\documentclass[12pt,reqno]{amsart}
\usepackage[english]{babel}
\usepackage[utf8]{inputenc}

\usepackage{amsmath, amsthm, amscd, amsfonts, amssymb, graphicx, color, mathrsfs}
\usepackage[bookmarksnumbered, colorlinks, plainpages]{hyperref}
\usepackage[all]{xy}
\usepackage[utf8]{inputenc}
\usepackage[english]{babel}
\usepackage{slashed}
\usepackage{cleveref}
\usepackage{mathtools}
\usepackage{amsrefs}
\usepackage{eqnarray}

\newtheorem{theorem}{Theorem}[section]
\newtheorem{lemma}[theorem]{Lemma}

\newtheorem{proposition}[theorem]{Proposition}
\newtheorem{corollary}[theorem]{Corollary}
\theoremstyle{definition}
\newtheorem{definition}[theorem]{Definition}

\theoremstyle{remark}
\newtheorem{remark}[theorem]{Remark}
\numberwithin{equation}{section}

\allowdisplaybreaks

\begin{document}
	\setcounter{page}{1}
	
	\title[]{On a class of pseudodifferential operators on the product of compact Lie groups}

	\author[S. Federico]{Serena Federico}
	\address{
		Serena Federico:
		\endgraf
		Department of Mathematics: Analysis Logic and Discrete
		\endgraf
		Mathematics
		\endgraf
		Ghent University
		\endgraf
		Krijgslaan 281, Ghent, B 9000
		\endgraf
		Belgium
		\endgraf
		{\it E-mail address} {\rm serena.federico@ugent.be}}
	
	\author[A. Parmeggiani]{Alberto Parmeggiani}
	\address{
		Alberto Parmeggiani:
		\endgraf
		Department of Mathematics
		\endgraf
		University of Bologna
		\endgraf
		Piazza di Porta S. Donato 5,
		40126 Bologna
		\endgraf
		Italy
		\endgraf
		{\it E-mail address} {\rm  alberto.parmeggiani@unibo.it}}
	\endgraf

	\thanks{Serena Federico has received funding from the European Unions Horizon 2020 research and innovation programme under the Marie Sk\l{}odowska-Curie grant agreement No 838661 and by the FWO Odysseus 1 grant G.0H94.18N: Analysis and Partial Differential Equations.}

	\begin{abstract} 
	In this paper a bisingular pseudodifferential calculus, along the lines of the one introduced by L.~Rodino in \cite{Rodino_ASNSP_1975}, is developed in the global setting of a product 
        of compact Lie groups. 
        The approach follows that introduced by M.~Ruzhansky and V.~Turunen \cite{MR-VT-Book_2010} (see also V.~Fischer \cite{Fischer-JFA_2015}), in that it exploits the harmonic analysis of the 
        groups involved.	
	\end{abstract}

\maketitle

\tableofcontents

\section{Introduction} \label{Intro}
In this paper we shall introduce a class of global pseudodifferential operators on the product of compact Lie groups and develope the corresponding global symbolic calculus in the spirit of the one introduced by Ruzhansky and Turunen in \cite{MR-VT-Book_2010} and of that introduced subsequently by Fischer in \cite{Fischer-JFA_2015}.

On the product of two manifolds the class we consider here was first studied by Rodino in \cite{Rodino_ASNSP_1975}, where, in particular, the author used the classical theory of pseudodifferential operators developed by H\"{o}rmander in \cite{Ho-FIO1} to construct an algebra of pseudodifferential operators containing the so-called bisingular operators.

The interest of our approach lies in the fact that it is {\it global} and based on the group structure and on its related representation theory.

As the classes $S^{m_1,m_2}(\Omega_1\times \Omega_2)$ in \cite{Rodino_ASNSP_1975} are not in general contained in any of the H\"{o}rmander classes $S^{m}(\Omega_1\times \Omega_2\times\mathbb{R}^{n_1+n_2})$, similarily in our case the classes $S^{m_1,m_2}(G_1\times G_2\times \widehat{G}_1\times\widehat{G}_2)$ are not in general contained in any class $S^m(G)$, with $G=G_1\times G_2$, defined by Ruzhansky and Turunen in \cite{MR-VT-Book_2010}. 

Natural examples of bisingular pseudodifferential operators in our setting, as in the general compact manifold setting, are tensor products of the form $A_1\otimes A_2$, where $A_i$ for $i=1,2$ is a pseudodifferential operator with symbol in the class $S^{m_i}(G_i)$ introduced in \cite{MR-VT-JW-JFAA_2014}, that is, $A_i\in L^m(G_i):=\mathrm{Op}(S^{m_i}(G_i))$, with $G_i$ being a compact Lie group.	

The study of these operators goes back to 1971, when Pilidi in \cite{Pilidi-DAN1971} reduced the boundary value problem for functions of two complex variables in bicylinders to the analysis of a bisingular equation on the two distinguished boundaries. In \cite{Pilidi-FA1971} the same author also developed a product calculus to deal with these objects and considered the corresponding index problem. 
Afterwards, a priori estimates and Fredholm properties for bisingular operators were studied by Rabinovi\u{c} in \cite{RAB-1973}, while in 1975 Rodino in \cite{Rodino_ASNSP_1975} introduced the so-called calculus of bisingular pseudodifferential operators.
Other related questions, such as residues and index problems, have been recently considered by Nicola and Rodino in \cite{MR-NR-2006}, while microlocal properties have been studied by Borsero and Schulz in \cite{BO-SH-2014}.

Let us also recall that a global version (i.e. in the Shubin setting of $\mathbb{R}^{n_1}\times\mathbb{R}^{n_2}$) of the calculus in \cite{Rodino_ASNSP_1975}
was developed by Battisti, Gramchev, Pilipovi\v{c} and Rodino in \cite{BA-GR-RO-PIL-2013}, and that other calculi of product type were developed by Dudu\v{c}ava in \cite{DUD-I-1979} 
and \cite{DUDI-II-1979}, and more recently by Melrose and Rochon in \cite{ME-ROC-2006}.

Note that a natural and immediate generalization of bisingular operators are the multisingular ones, whose prototype  are tensor products of the form $\bigotimes_{\substack{i=1}}^N A_i$, with $A_i\in L^{m_i}(G_i)$.  We will not pursue this topic here, but with suitable arrangements in the arguments used below one can define a multisingular pseudodifferential calculus on the direct product of 
finitely many compact Lie groups and define the corresponding multisingular pseudodifferential operators.

We want to remark that due to the intrinsic product stucture of the bisingular calculus, the suitable version of the celebrated G\aa{}rding inequality for elliptic operators is not available for {\it bielliptic} operators (see for instance \cite{RAB-1973}, where such inequality is attained only under very specific assumptions). Hence, it seems that for the class of bisingular operators, that serves as a model for degenerate elliptic operators, a more natural inequality to consider is the {\it sharp G\aa{}rding} inequality. We will analyze the problem of the validity of this inequality in a future paper that will be part II of the present work.

We finally conclude this introduction by giving the plan of the paper.

In Section 2 we shall recall some basic definitions on compact Lie groups, such as the notions of Fourier tranform, difference operators and Taylor expansion, as well as the standard quantization formula.

In Section 3 we introduce the class of bisingular symbols and define the corresponding pseudodifferential operators.

Section 4 will be devoted to the derivation of some fundamental kernel estimates needed to prove some asymptotic properties that are the object of Section 5.

Finally in Section 5 we develop the calculus, that is, we prove asymptotic formulas for the composition and for the adjoint of bisingular operators, and prove, after introducing ellipticity in the bisingular setting, the existence of parametrices for bielliptic operators.


\section{Preliminaries} 
In the sequel $G$ will be a compact Lie group, $\widehat{G}$ its unitary dual, that is the set of all equivalence classes of unitary representations of $G$, and $\mathrm{Rep}(G)$ the set of all the irreducible unitary representations of $G$. Since $G$ is compact, any given $\xi\in \mathrm{Rep}(G)$ is finite dimensional and we shall denote by $\mathcal{H}_\xi$ the associated representation space, and by $\mathcal{U}(\mathcal{H}_\xi)$ the corresponding space of unitary operators on $\mathcal{H}_\xi$.

The Fourier and inverse Fourier transforms on $G$ are given in terms of the representations of the group as follows.
\medskip

Given a function $f\in C^\infty(G)$, and $\xi\in \mathrm{Rep}(G)$, the (matrix-valued) global Fourier transform of $f$ at $\xi$ is defined by
$$\widehat{f}(\xi)=\int_G f(x)\xi^*(x)dx,$$
where $\xi^*(x):=\,^t\overline{\xi(x)}$ stands for the adjoint representation of $\xi$, while $dx$ denotes the Haar measure on the group. Notice that, given $\xi: \mathcal{H}_\xi\rightarrow \mathcal{U}(\mathcal{H}_\xi)$, and $d_\xi:=\mathrm{dim}(\xi):=\mathrm{dim}(\mathcal{H}_\xi)$, then  $\widehat{f}(\xi)\in\mathbb{C}^{d_\xi\times d_\xi}$.
Correspondingly, the  inverse Fourier transform is given by 
$$ f(x)=\sum_{[\xi]\in\widehat{G}}d_\xi \mathrm{Tr}(\xi(x)\widehat{f}(\xi)),$$
where $\mathrm{Tr}(A)$ denotes the trace of the matrix $A$.  

Related to the previous formulas one has the following Parseval identity 
$$\| f\|^2_{L^2(G)}=\sum_{[\xi]\in\widehat{G}} d_\xi \|\widehat{f}(\xi)\|^2_{HS}=: \|\widehat{f}\|_{\ell^2(\widehat{G})}^2,$$
where $\|\widehat{f}(\xi) \|_{HS}:=\left(\mathrm{Tr}(\widehat{f}(\xi)\widehat{f}(\xi)^*)\right)^{1/2}$ is the Hilbert-Schmidt norm .

\medskip

In order to deal with (matrix-valued) functions on $\widehat{G}$ we will need to make use of the so called difference operators that we next define following \cite{MR-VT-JW-JFAA_2014}.

\begin{definition}\label{diff_1}
We say that $Q_\xi$ is a difference operator of order $k$ on $\mathcal{F}(\mathcal{D}'(G))$ (the image of the group Fourier transform of distributions on $G$) if
$$Q_\xi \widehat{f}(\xi)=\widehat{q_Qf}(x),$$
for a function $q_Q\in C^\infty(G)$ vanishing of order $k$ at the identity element $e$ of $G$, that is, $q_Q$ is  such that $q_Q(e)=P_x q_Q(e)=0$ for all left-invariant differential operators $P_x\in \mathrm{Diff^{k-1}(G)}$ of order $k-1$. 
\end{definition}

We shall denote by $\mathrm{diff}^k(\widehat{G})$ the set of all difference operators of order $k$ on $\widehat{G}$.

\begin{definition}\label{admissible-coll}
	A collection of $n_\triangle\geq n=\mathrm{dim}(G)$  difference operators $\triangle_1,...,\triangle_{n_\triangle}$ in $\mathrm{diff}^1(\widehat{G})$ is called {\it admissible} if the corresponding functions $q_1,...,q_{n_\triangle} \in C^\infty(G)$ are such that $q_1(e)=\ldots =q_{n_\triangle}(e)=0$, and  $dq_j(e)\neq 0$ for all $j=1,...,{n_\triangle},$ with $\mathrm{rank}(dq_1(e),...,dq_{n_\triangle}(e))=n$. Finally, a collection of difference operators is called {\it strongly admissible} if $\bigcap_j\{x\in G; q_j(x)=0\}=\{e\}$.
\end{definition}

Given a fixed family of functions $Q=\{q_j\}_{j=1,...,{n_\triangle}}$, we shall denote by 
\begin{itemize}
\item  $\triangle_Q$ the associated admissible collection of difference operators;
\item  $q^\alpha:=q_1^{\alpha_1}\ldots q_{n_\triangle}^{\alpha_{n_\triangle}}$;
\item  $\triangle_{Q,j}=\triangle_{q_j}$ and $\triangle_Q^\alpha:=\triangle_{Q,1}^{\alpha_1}\cdots \triangle_{Q,n_\triangle}^{\alpha_{n_\triangle}}$ the corresponding element in $\mathrm{diff}^{|\alpha|}(\widehat{G})$.
\end{itemize}
Additionally, once the collection of difference operators is fixed, namely the corresponding family of functions $Q$ is fixed, one can find a family of differential operators in $\mathrm{Diff}^{|\alpha|}(G)$, denoted by $\partial_x^{(\alpha)}$, such that the following Taylor's formula holds
$$f(x)=\sum_{|\alpha|<N}\frac{1}{\alpha!}q(x)^\alpha\partial_x^{(\alpha)}f(e)+\mathcal{O}(h(x)^N),\quad h(x)\rightarrow 0,$$
for all $f\in C^\infty(G)$, where $h(x)$ is the geodesic distance from $x$ to $e_G$.
The differential operators $\partial_x^{(\alpha)}$ can be replaced by $\partial^\alpha_x:=\partial_{x_1}^{\alpha_1}\ldots \partial_{x_n}^{\alpha_n}$, with $\partial_{x_j}$, $j=1,\ldots, n$, being a collection of left-invariant first order differental operators corresponding to some linearly independent left-invariant vector fields on $G$ ($\partial_{x_j}$ are not the Euclidean directional derivatives here).

\begin{remark}\label{rmk.vf}
Note that we are assuming the Lie algebra $\mathfrak{g}$ to be the space of left-invariant vector fields. In particular, we shall use the notation $\partial_{x_j}$ and $\tilde{\partial}_{x_j}$ for the left and right invariant vector fields, respectively.
Once we fix an orthonormal basis of left-invariant vector fields for $\mathfrak{g}$, then any element of $\mathrm{Diff^{k}(G)}$ (the space of smooth vector fields on $G$) can be written as a linear combination in terms of the elements of the basis. Note also that a similar property holds for right-invariant vector fields. 
\end{remark}

By Lemma 4.4 in \cite{MR-VT-JW-JFAA_2014}, the family of functions $\{q_{ij}=\xi_{ij}-\delta_{ij}\}_{[\xi]\in \widehat{G},1\leq i,j\leq d_\xi}$ always induces a strongly admissible collection of difference operators, therefore we choose the latter as the fixed admissible collection for the rest of the paper. 

In the context of the difference operators defined above, the following notion of Leibniz formula is adopted (see \cite{Fischer-JFA_2015})

\begin{definition}\label{Leibniz-like formula}
A collection of $\triangle=\triangle_{Q}$ difference operators satisfies the {\it Leibniz-like} property if, for any Fourier transforms $\widehat{f_1}$ and  $\widehat{f_2}$ (with $f_1,f_2\in\mathcal{D}'(G)$),
$$\triangle_{Q,j}(\widehat{f_1}\widehat{f_2})=\triangle_{Q,j}(\widehat{f_1})\widehat{f_2}+
\widehat{f_1}\triangle_{Q,j}(\widehat{f_2})+\sum_{1\leq l,k\leq n_\triangle}c_{l,k}^{(j)}\triangle_{Q,l}(\widehat{f_1})\triangle_{Q,k}(\widehat{f_2}),$$
for some coefficients $c_{l,k}^{(j)}\in \mathbb{C}$ depending only on $l,k,j$ and $\triangle$.
\end{definition}
If $\triangle$ is a collection satisfying the Leibniz-like formula, then, recursively, for any given $\alpha\in \mathbb{N}^{n_\triangle}_0$, one has
\begin{equation}\label{Leinbniz}
\triangle_{Q}^\alpha(\widehat{f_1}\widehat{f_2})=\sum_{|\alpha|\leq |\alpha_1|+|\alpha_2|\leq 2|\alpha|} c_{\alpha_1,\alpha_2}^\alpha (\triangle_{Q}^{\alpha_1}\widehat{f_1})(\triangle_{Q}^{\alpha_2}\widehat{f_2}),
\end{equation}
for some coefficients $c_{\alpha_1,\alpha_2}^\alpha\in \mathbb{C}$ depending on $\alpha_1,\alpha_2,\alpha$ and $\triangle$, with $c^\alpha_{\alpha,0}=c^\alpha_{0,\alpha,}=1$.

\begin{remark}\label{rmk-Leibniz}
We remark that there always exists a strongly admissible collection of difference operators as in Definition \ref{admissible-coll} for which the Leibniz property above holds (see Corollary 5.13 in \cite{Fischer-JFA_2015}). In particular this is always the case for the strongly admissible collection $\triangle_{Q}$ with $Q=\{q_{ij}=\xi_{ij}-\delta_{ij}, 1\leq i,j\leq d_\xi, \xi\in \widehat{G}\}$.
\end{remark}
\medskip 

In order to introduce the precise difference operators we are going to exploit in our setting, it is more convenient for us to use the definition of difference operators given in \cite{Fischer-JFA_2015}, which, in turn, gives rise to the same difference operators defined above by means of the admissible collection $\{q_{ij}=\xi_{ij}-\delta_{ij}\}_{[\xi]\in \widehat{G},1\leq i,j\leq d_\xi}$.

\begin{definition}\label{diff_2}
For each $\tau,\xi\in\mathrm{Rep}(G)$ we define the linear mapping $\triangle_\tau \widehat{f}(\xi)$ on $\mathcal{H}_\tau\otimes\mathcal{H}_\xi$ by
\begin{equation}\label{diff_2-formula}
	\triangle_\tau \widehat{f}(\xi)=\widehat{f}(\tau\otimes\xi)-\widehat{f}(I_{d_\tau}\otimes \xi).
\end{equation}
The restriction of $\triangle_\tau\widehat{f}(\xi)$ to any occurrence of $\rho\in\widehat{G}$ in the decomposition into irreducibles of $\tau\otimes\xi$ defines the same mapping on $\mathcal{H}_\rho$, while the restriction to any $\rho\in\widehat{G}$ not appearing in the decomposition of $\tau\otimes\xi$ is fixed to be zero. With these conventions the operation $\triangle_\tau$ is called {\it difference operator} associated with $\tau\in\mathrm{Rep}(G)$.
\end{definition}

\begin{definition}\label{def.admissible.relative}
Let $G=G_1\times G_2$ a compact Lie group such that $G_i$ is compact for all $i=1,2$, and $n_i:=\mathrm{dim}(G_i)$. Let $e=(e_1,e_2)$ be
the neutral element of $G$. A collection of $n_P:=n_{\triangle_P}\geq n_1$ of difference operators $\triangle_{p_1},\ldots,\triangle_{p_{n_{P}}}\in \mathrm{diff}^1(\widehat{G})$ is called {\it admissible} relative to $G_1$ if the corresponding functions $p_1,\ldots,p_{n_{P}}\in C^\infty(G)$ are such that $p_1(e)=\ldots=p_{n_P}(e)=0$, and $dp_j(e)\neq 0$ for all $j=1,...,n_{P}$, with $\mathrm{rank}(dp_1(e),\ldots,dp_{n_{P}}(e))=n_1$. 

The collection is called {\it strongly admissible} relative to $G_1$ if $\bigcap_{j=1}^{n_{\triangle_P}}\{ x\in G; p_j(x)=0\}=\{e_1\}\times G_2.$
Admissible and strongly admissible collections relative to $G_2$ are defined similarily by reversing the role of $G_1$ and $G_2$.
\end{definition}

We then consider the family of functions 
\begin{equation}\label{r}
R=\{r_{ij}^{(\tau)};1\leq i,j\leq d_\tau,  \tau\in \widehat{G}\}=\{(I_{d_{\tau_1}}\otimes\tau_2-I_{d_\tau})_{ij}; 1\leq i,j\leq d_\tau,  \tau\in \widehat{G}\} 
\end{equation}
and 
\begin{equation}\label{p}
P=\{p_{ij}^{(\tau)}; 1\leq i,j\leq d_\tau,  \tau\in \widehat{G}\}=\{(\tau_1\otimes I_{d_{\tau_2}}-I_{d_\tau})_{ij}; 1\leq i,j\leq d_\tau,  \tau\in \widehat{G}\},
\end{equation}
so that both $\triangle_R$ and $\triangle_P$ are strongly admissible collections relative to $G_1\simeq G_1\times \{e_2\}\subset G$ and $G_2\simeq \{e_1\}\times  G_2\subset G$ respectively (see Definition \ref{def.admissible.relative}). After ordering the (huge but finite) families $P$ and $R$ above, that is, writing
$$P=\{p_k, k=1,...,n_{{P}} \}, \quad R=\{r_k, k=1,...,n_{{R}} \}, $$
where each $p_k, r_k$  are functions of the form $p_{ij}^{(\tau)}, r_{ij}^{(\tau)}$ respectively, for some $\tau\in\widehat{G}$ and some $i,j\in \{0,...,d_\tau\}$,
we may define
\begin{equation}\label{difference-pr}
\triangle^{\alpha,\beta}:=\triangle_P^\alpha \triangle_R^\beta =\triangle^{\alpha_1}_{p_1}\cdots\triangle^{\alpha_{n_{P}}}_{p_{n_{P}}}\,
\triangle^{\beta_1}_{r_1}\cdots\triangle^{\beta_{n_{R}}}_{r_{n_{R}}}.
\end{equation}
Note that the function $p_j$, for all $j=1,...,n_{{P}},$ is independent of $x_2\in G_2$, and, similarily, the function $r_j$, for all $j=1,...,n_{{R}}$, is independent of $x_1$.

These difference operators, namely of the form \eqref{difference-pr}, will be the ones used for the rest of the paper. Note that $\triangle_{P}^\alpha$ and $\triangle_{Q}^\beta$
may be tought of as ''partial difference operators'' in the ``directions'' of $\widehat{G}_1$ and $\widehat{G}_2$ respectively.

By Remark \ref{rmk-Leibniz} and formula \eqref{Leibniz-like formula} applied iteratively, we have the following Leibniz-like formula for the difference operators we are considering.

\begin{proposition}\label{prop-product-Leibniz-formula}
Let $G=G_1\times G_2$ , with $G_1$ and $G_2$ compact Lie groups. Then, for any $\alpha\in \mathbb{N}_0^{n_{\triangle_{P}}}$, $\beta\in \mathbb{N}_0^{n_{\triangle_{R}}}$, and for all Fourier transforms $\widehat{f_1}, \widehat{f_2}$ (with $f_1,f_2\in\mathcal{D}'(G)$), we have
\begin{equation}\label{product-Leibniz-formula}
\triangle^{\alpha,\beta}(\widehat{f_1}\widehat{f_2})=\sum_{|\alpha|\leq |\alpha_1|+|\alpha_2|\leq 2|\alpha|}\sum_{|\beta|\leq |\beta_1|+|\beta_2|\leq 2|\beta|} c^\alpha_{\alpha_1,\alpha_2}
c^{\beta}_{\beta_1,\beta_2} (\triangle^{\alpha_1,\beta_1}\widehat{f_1})(\triangle^{\alpha_2,\beta_2}\widehat{f_2}),
\end{equation}
for some coefficients $c^\alpha_{\alpha_1,\alpha_2}, c^{\beta}_{\beta_1,\beta_2}\in\mathbb{C}$ such that $c^{\beta}_{\beta,0}=c^{\beta}_{0,\beta}= c^\alpha_{\alpha,0}= c^\alpha_{0,\alpha}=1$.
\end{proposition}
\proof
The proof immediately follows by application of formula \eqref{Leinbniz} twice, that is, for $\triangle^\alpha_{P}$ and $\triangle_{Q}^\beta$ respectively.
\endproof

Observe now that, since the families of functions in \eqref{r} and \eqref{p} defining admissible collections of difference operators on $\widehat{G}$ relative to $G_1$ and $G_2$ are fixed, on denoting by $x=(x_1,x_2)$ an element of $G=G_1\times G_2$, with $\mathrm{dim}(G_1)=n_1$ and $\mathrm{dim}(G_2)=n_2$, we can find a family of differential operators
$$\partial^{\alpha,\beta}_x:=\partial_{x_1}^{\alpha}\partial^\beta_{x_2}$$ such that the following form of Taylor's formula holds (see, for instance, \cite{Fischer-JFA_2015})
$$f(x)=\sum_{|\alpha|<N}\sum_{|\beta|<N}\frac{1}{\alpha!\beta!}q^{\alpha,\beta}(x^{-1})\partial_x^{\alpha,\beta}f(e)+\sum_{\substack{|\alpha+\beta|=2N, \\
|\alpha|\geq N \vee |\beta|\geq N}}\frac{1}{\alpha!\beta!}q^{\alpha,\beta}(x^{-1})f_{\alpha,\beta}(x),$$
where 
$$q^{\alpha,\beta}(x):=r(x)^\alpha p(x)^\beta=r_{1}(x)^{\alpha_1}\ldots r_{n_{R}}(x)^{\alpha_{n_{R}}}\,
p_{1}(x)^{\beta_1}\ldots p_{n_{P}}(x)^{\beta_{n_{P}}}.$$
Recall that, in particular, we will have 
${\triangle_{R,j}}\widehat{f}(\xi):=\widehat{r_{j}f}(\xi)$ and ${\triangle_{P,k}}=\widehat{p_{k}f}(\xi)$. Moreover, the differential operators $\partial^{\alpha,\beta}_x$ are chosen so as to satisfy
$\partial_{x_1}^\alpha(p(x)^\alpha)=\partial_{x_2}^\beta(r(x)^\beta)=1$ for all $\alpha,\beta$ such that $|\alpha|=|\beta|=1$. In particular, since $P$ and $R$ are strongly admissible collections relative to $G_1$ and $G_2$ respectively, we have that there are $n_1$ and $n_2$ elements in $P$ and $R$ respectively,  
say $(p_1,\ldots, p_{n_1})$ and $(r_1,\ldots, r_{n_2})$, such that $(\partial_{x_{1,1}},\ldots, \partial_{x_{1,n_1}}, \partial_{x_{2,1}},\ldots, \partial_{x_{2,n_2}})$ can be identified with $(dp_1(e),\ldots, dp_{n_1}(e), dr_1(e),\ldots, dr_{n_2}(e))$ (where $df(e)$ denotes the differential computed at $e$) by duality, and we fix the former as the basis of the Lie algebra $\mathfrak{g}$. 
We stress that the choice of $q^{\alpha,\beta}(x^{-1})$ istead of $q^{\alpha,\beta}(x)$ is technical (see \cite{MR-VT-Book_2010}). Note finally that the formula above can be derived by application of Taylor's formula twice, that is, first with respect to the variable $x_1$ by using the functions $q^{\alpha,0}(x)=q^{\alpha,0}(x_1)$, and then by expanding again with respect to $x_2$ and using $q^{0,\beta}(x)=q^{0,\beta}(x_2)$.

\section{Bisingular symbols on $G=G_1\times G_2$}
In this section we define what we shall call class of {\it bisingular symbols}, since, as pointed out by L.~Rodino in \cite{Rodino_ASNSP_1975}, it contains symbols of operators of bisingular type (see \cite{Pilidi-DAN1971},\cite{Pilidi-FA1971} and \cite{Simonenko-FAP-1971}).
\medskip

\noindent\textbf{Notation.} In what follows we call $x=(x_1,x_2)$ an element of $G=G_1\times G_2$ and $\xi:=\xi_1\otimes \xi_2$ an element of $\widehat{G}$, where $\xi_j\in \widehat{G}_j$. By using the definitions above, and fixing the families $R$ and $P$, we define $\triangle_1^\alpha:=\triangle_{P}^\alpha$, $\triangle_{2}^\beta:=\triangle_{Q}^\beta$, and $\partial_1^\alpha:=\partial^\alpha_{x_1}=\partial^{\alpha_1}_{x_{11}}\ldots \partial^{\alpha_{n_1}}_{x_{1n_1}}$, $\partial_2^\beta:=\partial_{x_2}^\beta=\partial^{\beta_{1}}_{x_{21}}\ldots \partial^{\beta_{n_2}}_{x_{2n_2}}$ as above (where, as previously mentioned, $\partial_{x_j}$ \textit{are not} the Euclidean directional derivatives). We shall finally put $\partial^{\alpha,\beta}:=\partial^\alpha_{x_1}\partial^\beta_{x_2}$  and analogously for the difference operators $\triangle^{\alpha,\beta}$.
\medskip

We recall that, given a continuous linear operator $A$ from $C^\infty(G)$ to $\mathcal{D}'(G)$, its matrix-valued symbol $\sigma_A(x,\xi)\in \mathbb{C}^{d_\xi\times d_\xi}$ (as introduced in \cite{MR-VT-Book_2010}) is given by
\begin{equation}\label{symbol}
 \sigma_A(x,\xi)=\xi^*(x) (A\xi)(x),
\end{equation}
and that
$$Af(x)=\sum_{[\xi]\in \widehat{G}} d_\xi \mathrm{Tr}(\xi(x)\sigma_A(x,\xi)\widehat{f}(\xi)), \quad f\in C^\infty(G),$$
holds in the sense of distributions and the sum is independent of the choice of the rapresentative $\xi$ of the class $[\xi]$.
 
\begin{definition}
Let $G=G_1\times G_2$ be a compact Lie group and define $n_i:=\mathrm{dim}(G_i)$. We call \textit{class of bisingular symbols of order} $(m_1,m_2)\in \mathbb{R}^2$ the set $S^{m_1,m_2}(G\times \widehat{G})$ of all $a: G\times \widehat{G}\longrightarrow \bigcup_{[\xi]\in \widehat{G}} \mathbb{C}^{d_\xi\times d_\xi}$  that are smooth in $x\in G$ and such that, for all multiindices $\alpha_1\in\mathbb{N}_0^{n_1},\alpha_2\in\mathbb{N}_0^{n_2}, \beta_1\in\mathbb{N}_0^{\triangle_P},\beta_2\in\mathbb{N}_0^{\triangle_R}$, 
$$\| \partial^{\alpha_1}_{x_1}\partial^{\alpha_2}_{x_2}\triangle_1^{\beta_1}\triangle_2^{\beta_2}
a(x_1,x_2,\xi_1,\xi_2)\|_{op}\leq C_{\alpha_1,\alpha_2,\beta_1,\beta_2}\langle\xi_1\rangle^{m_1-|\beta_1|}\langle\xi_2\rangle^{m_2-|\beta_2|},
$$
where
$$ \| a\|_{op}:=\sup\{ | a(x,\xi)v|_{\ell^2}; v\in \mathbb{C}^{d_\xi}, |v|_{\ell^2}\leq 1\}.$$
Additionally we shall denote by $S^{-\infty,-\infty}(G\times  \widehat{G}):=\bigcap_{(m_1,m_2)\in \mathbb{Z}^2} S^{m_1,m_2} (G\times  \widehat{G})$ the class of smoothing elements.
\end{definition}

It is important to bear in mind that $\widehat{G}=\widehat{G}_1\times\widehat{G}_2$.

Due to the equivalence of $\|a \|_{\mathscr{L}(\mathcal{H}_\xi)}$ and $\|a\|_{op}$, we will freely use both notations below.

Let us remark that, as in the standard case, the space $S^{m_1,m_2}(G\times\widehat{G})$ is a Fréchet space equipped with the seminorms
$$\|\sigma\|_{S^{m_1,m_2}_{(a_1,a_2),(b_1,b_2)}}:= \max_{\substack{|\alpha_1|\leq a_1, |\alpha_2|\leq a_2\\  |\beta_1|\leq b_1,|\beta_2|\leq b_2}}\sup_{(x,\xi)\in G\times\widehat{G}} \langle\xi_1\rangle^{-m_1+|\alpha_1|} \langle\xi_2\rangle^{-m_2+|\alpha_2|} \| \triangle^{\alpha_1,\alpha_2}\partial_{x_1}^{\beta_1,\beta_2} \sigma (x,\xi)\|_{\mathscr{L}(\mathcal{H}_\xi)},$$
with $a_1,a_2,b_1,b_2\in\mathbb{N}_0$.

To each matrix-valued symbol $a\in S^{m_1,m_2}(G\times \widehat{G}_1\times \widehat{G}_2)$ one can associate an operator $\mathrm{Op}(a)$ by means of the following quantization formula
$$\mathrm{Op}(a)\varphi(x):=\sum_{[\xi]\in \widehat{G}}d_\xi\mathrm{Tr}(\xi(x)a(x,\xi)\widehat{\varphi}(\xi))$$
\begin{equation}\label{quant-formula}
	=\sum_{[\xi_1]\in \widehat{G}_1}\sum_{[\xi_2]\in \widehat{G}_2}d_{\xi_1}d_{\xi_2}\mathrm{Tr}((\xi_1\otimes\xi_2)(x)a(x,\xi_1,\xi_2)
	\widehat{\varphi}(\xi_1\otimes \xi_2)),
	\end{equation}
and we shall denote by $L^{m_1,m_2}(G)$ the class of operators of the previous form, that is, those obtained by quantizing symbols in $S^{m_1,m_2}(G\times \widehat{G})$ as in \eqref{quant-formula}. These operators will be called  {\it binsingular operators} of order $(m_1,m_2)$ on $G=G_1\times G_2$.

Moreover, with any $a\in S^{m_1,m_2}(G\times \widehat{G})$, we associate the maps

$$ G_1\times \widehat{G}_1\ni (x_1,\xi_1)\longmapsto a(x_1,x_2,\xi_1,D_2)\in L^{m_2}(G_2),$$
$$G_2\times \widehat{G}_2\ni (x_2,\xi_2)\longmapsto a(x_1,x_2,D_1,\xi_2)\in L^{m_1}(G_1),$$
where $L^{m_1}(G_1)$ and $L^{m_2}(G)$ are classes of operators on $G_1$ and $G_2$ respectively obtained by means of the quantization formulas

$$a(x_1,x_2, \xi_1, D_2)\varphi(x_2)=\sum_{[\xi_2]\in \widehat{G}_2}d_{\xi_2}\mathrm{Tr}\left( (I_{d_{\xi_1}}\otimes\xi_2(x_2))a(x_1,x_2,\xi_1,\xi_2)
 \times  (I_{d_{\xi_1}}\otimes \widehat{\varphi }(\xi_2)) \right)$$
and
$$a(x_1,x_2, D_1, \xi_2)\varphi(x_1)=\sum_{[\xi_1]\in \widehat{G}_1}d_{\xi_1}\mathrm{Tr}\left( (\xi_1(x_1)\otimes I_{d_{\xi_2}})a(x_1,x_2,\xi_1,\xi_2) 
 \times ( \widehat{\varphi }(\xi_1) \otimes I_{d_{\xi_2}}) \right).$$

It  is important to stress that  the symbol $a\in S^{m_1,m_2}(G\times \widehat{G})$ is uniquely determined by one of these maps.

Troughout the paper we will often write $a(x,\xi)$ in place of $a(x_1,x_2,\xi_1,\xi_2)$, where $\xi=\xi_1\otimes\xi_2\in \widehat{G}$, and $S^{m_1,m_2}(G\times \widehat{G}_1\times \widehat{G}_2)$  in place of $S^{m_1,m_2}(G\times \widehat{G})$.

\begin{remark}
Notice that, in general, there is no $m\in \mathbb{R}$ such that $S^{m_1,m_2}(G\times \widehat{G})\subset S^m(G\times \widehat{G})$. However we always have that $S^{m_1,m_2}(G\times \widehat{G})\subset S^m_{0,0}(G\times \widehat{G})$ for some $m\in \mathbb{R}$.
\end{remark}

Given a continuous linear operator $A:\mathcal{D}(G)\rightarrow\mathcal{D}'(G)$, (where $\mathcal{D}(G):=C^\infty(G)$), its right-convolution kernel $R_A\in \mathcal{D}'(G\times G)$ is defined by
\begin{equation}\label{Kernel}
A\varphi(x)=\int_G \varphi(y) R_A(x,y^{-1}x) dy=(R_A(x,\cdot)\ast \varphi )(x).
\end{equation}
Therefore, given $A\in L^{m_1,m_2}(G)$ with simbol $\sigma_A\in S^{m_1,m_2}(G\times\widehat{G})$, one has 
$$\sigma_A(x,\xi):=(\mathcal{F}_{y\rightarrow\xi}\,R_A)(x,\xi),$$
where
$$R_A (x,y):=\sum_{[\xi]\in \widehat{G}}
d_{\xi}\mathrm{Tr}\left(\xi(y)a(x,\xi) \right),$$
with $y=(y_1,y_2)\in G=G_1\times G_2$ and $\xi\in\widehat{G}$ of the form $\xi=\xi_1\otimes\xi_2$, with $(\xi_1,\xi_2)\in \widehat{G}_1\times\widehat{G}_2$.

For any fixed $(x_1,\xi_1)\in G_1\times \widehat{G}_1$ and $(x_2,\xi_2)\in G_2\times \widehat{G}_2$, we can write, respectively, the operators $a(x_1,x_2,\xi_1,D_2)$ and $a(x_1,x_2,D_1,\xi_2)$ defined above in terms of their (right-)convolution kernels, that is,
$$a(x_1,x_2,\xi_1,D_2)\varphi(x_2)=(R^2_a(x_1,x_2,\xi_1, \cdot)\ast_{G_2}\varphi)(x_2)$$
$$a(x_1,x_2,D_1,\xi_2)\varphi(x_1)=(R^1_a(x_1,x_2,\cdot, \xi_2)\ast_{G_1}\varphi)(x_1)$$
where
\begin{equation}\label{kernel1}
R_a^2 (x_1,x_2,\xi_1, y_2):=\sum_{[\xi_2]\in\widehat{G}_2}d_{\xi_2}\mathrm{Tr}
\left( (I_{\xi_1}\otimes\xi_2(y_2)) a(x_1,x_2, \xi_1,\xi_2)\right)
\end{equation}
and
\begin{equation}\label{kernel2}
R_a^1 (x_1,x_2,y_1, \xi_2):=\sum_{[\xi_1]\in\widehat{G}_1}d_{\xi_1}\mathrm{Tr}
\left( (\xi_1(x_1)\otimes I_{\xi_2}) a(x_1,x_2, \xi_1,\xi_2)\right).
\end{equation}
Due to the orthogonality property of irriducible representations we have that
$$a(x_1, x_2,\xi_1,\xi_2)=\int_{G_1}  R_a^1 (x_1,x_2,y_1, \xi_2) (\xi_1(y_1)^*\otimes I_{\xi_2}) d y_1$$
and
$$a(x_1, x_2,\xi_1,\xi_2) =\int_{G_2} R_a^2 (x_1,x_2,\xi_1, y_2) (I_{\xi_1}\otimes\xi_2(y_2)^*)  dy_2.$$

\begin{definition}\label{partial_composition}
Given $a\in S^{m_1,m_2}(G\times \widehat{G})$ and $b\in S^{m'_1,m'_2}(G\times \widehat{G})$, we shall denote by $(a\circ_{\xi_1}b)(x_1,x_2,\xi_1, \xi_2)$ and $ (a\circ_{\xi_2}b)(x_1,x_2,\xi_1, \xi_2)$ the symbols in $S^{m_1+m'_1,m_2+m'_2}(G\times \widehat{G})$  corresponding to the operators

$$ (a\circ_{\xi_1}b)(x_1,x_2,D_1,\xi_2)\varphi(x_1)=
a(x_1,x_2,D_1,\xi_2)b(x_1,x_2,D_1,\xi_2)\varphi(x_1),\quad \forall \varphi \in C^\infty(G_1),$$
  and
 $$ (a\circ_{\xi_2}b)(x_1,x_2,\xi_1,D_2)\psi(x_2)=
 a(x_1,x_2,\xi_1,D_2)b(x_1,x_2,\xi_1,D_2)\psi(x_2),\quad \forall \psi \in C^\infty(G_2).$$
\end{definition}

By considering the right-convolution kernels it is not difficult to show that 
$$(a\circ_{\xi_1}b)(x_1,x_2,\xi_1,\xi_2)\sim \sum_{ |\alpha_1|\geq 0} (\triangle^{\alpha_1,0}a(x,\xi))\,\partial^{\alpha_1,0}b(x,\xi)$$
and
$$(a\circ_{\xi_2}b)(x_1,x_2,\xi_1,\xi_2)\sim \sum_{ |\alpha_2|\geq 0} (\triangle^{0,\alpha_2}a(x,\xi))\,\partial^{0,\alpha_2}b(x,\xi),$$
where for all $N>0$ we have
$$r_N^1(x,\xi):=(a\circ_{\xi_1}b)(x_1,x_2,\xi_1,\xi_2)-\sum_{ |\alpha_1|<N} (\triangle^{\alpha_1,0}a(x,\xi))\,\partial^{\alpha_1,0}b(x,\xi)=\sum_{ |\alpha_1|=N} (\triangle^{\alpha_1,0}a(x,\xi)) b_{\alpha_1}(x,\xi)$$

$$r_N^2(x,\xi):=(a\circ_{\xi_2}b)(x_1,x_2,\xi_1,\xi_2)-\sum_{ |\alpha_2|\geq 0} (\triangle^{0,\alpha_2}a(x,\xi))\,\partial^{0,\alpha_2}b(x,\xi)=\sum_{ |\alpha_2|=N} (\triangle^{0,\alpha_2}a(x,\xi) )b_{\alpha_2}(x,\xi),$$
for suitable $ b_{\alpha_1}, b_{\alpha_2}$ having the same properties as $b$, that is, $b_{\alpha_1}, b_{\alpha_2}\in S^{m'_1,m'_2}(G\times \widehat{G})$.

\medskip

Let $a\in S^{m_1,m_2}(G\times \widehat{G})$ and denote by $\mathrm{Op}(a_{(x_2,\xi_2)})(x_1,D_1):=a(x_1,x_2,D_1,\xi_2)$
the operator defined above and belonging to $L^{m_1}(G_1)$ for all $(x_2,\xi_2)\in G_2\times \widehat{G}_2$. Then it is possible to define the adjoint of $\mathrm{Op}(a_{(x_2,\xi_2)})$ (as an operator on $G_1$), denoted by $\mathrm{Op}(a_{(x_2,\xi_2)})(x_1,D_1)^{*_1}:=a(x_1,x_2,D_1,\xi_2)^{*_1}$, as the operator satisfying
\begin{equation}\label{star1}
(\mathrm{Op}(a_{(x_2,\xi_2)})u,v)_{L^2(G_1)}=(u,\mathrm{Op}(a_{(x_2,\xi_2)})^{*_1}v)_{L^2(G_1)},\quad u,v\in\mathcal{D}(G_1),
\end{equation}
where $(\cdot,\cdot)_{L^2(G_1)}$ stands for the scalar product on $L^2(G_1)$.

In a similar way, on denoting by
 $$\mathrm{Op}(a_{(x_1,\xi_1)})(x_2,D_2):=a(x_1,x_2,\xi_1,D_2)$$  the operator belonging to $L^{m_2}(G_2)$ for all $(x_1,\xi_1)\in G_1\times \widehat{G}_1$, one can define the adjoint operator $\mathrm{Op}(a_{(x_1,\xi_1)})(x_2,D_2)^{*_2}:=a(x_1,x_2,\xi_1,D_2)^{*_2}$
as the one satisfying 
\begin{equation}\label{star2}
(\mathrm{Op}(a_{(x_1,\xi_1)})u,v)_{L^2(G_2)}=(u,\mathrm{Op}(a_{(x_1,\xi_1)})^{*_2}v)_{L^2(G_2)},\quad u,v\in\mathcal{D}(G_2),
\end{equation}
with $(\cdot,\cdot)_{L^2(G_2)}$ denoting the scalar product on $L^2(G_2)$.
\medskip

\noindent \textbf{Sobolev spaces $H^{s_1,s_2}(G)$.}
We shall now define what we shall call {\it bisingular Sobolev spaces} which are the ones to be naturally 
used in this setting. We will also see that bisingular operators exhibit continuity properties when acting on these spaces.

Let us consider the operator $L$ on $G=G_1\times G_2$, defined as
$$L:= (I_1+L_{G_1})\otimes(I_2+L_{G_2}),$$
where $L_{G_i}$ and $I_i$ denote the positive Laplace operator and the identity operator on  $G_i$ respectively. 

The operator $L$ will be called {\it bilaplacian}, since, as expected, it will play the role of the Laplacian in this setting.

By formula \eqref{symbol} we have that the symbol of the bilaplacian is given by
$$\sigma_L(\xi)=\sigma_L(\xi_1\otimes\xi_2)=\langle \xi_1\rangle^2\langle \xi_2\rangle^2 I_{d_\xi},$$
where $\langle \xi_i\rangle:=(1+\lambda^2_{\xi_1})^{1/2}$, with $\lambda^2_{\xi_i}>0$ being the eigenvalue of $L_{G_i}$ relative to the representation 
$\xi_i\in\widehat{G}_i$, and $I_{d_{\xi_i}}\in \mathbb{C}^{d_{\xi_i}\times d_{\xi_i}}$ is the identity matrix.

\begin{definition}[Bisingular Sobolev space of order $(s_1,s_2)$] \label{def-bisingular Sobolev space}
We shall call bisingular Sobolev space of order $(s_1,s_2)$ the space
$$H^{s_1,s_2}(G):=\{f\in \mathcal{D}'(G); \langle \xi_1\rangle^{s_1}  \langle \xi_2\rangle^{s_2}\widehat{f}(\xi)\in \ell^2(\widehat{G})\},$$
equipped with the norm
$$\| f\|_{s_1,s_2}:=\left(\sum_{[\xi]\in \widehat{G}} d_\xi \langle \xi_1\rangle^{2s_1}  \langle \xi_2\rangle^{2s_2} \mathrm{Tr}(\widehat{f}(\xi)^* \widehat{f}(\xi)  )\right)^{1/2}$$
$$= \|  \langle \xi_1\rangle^{s_1}  \langle \xi_2\rangle^{s_2}\widehat{f}\|_{\ell^2(\widehat{G})}=: \| \widehat{f}\|_{h^{s_1,s_2}(\widehat{G})},$$
where 
$$h^{s_1,s_2}(\widehat{G}):=\{ \widehat{f}\in \mathcal{F}(\mathcal{D}'(G));  \langle \xi_1\rangle^{s_1}  \langle \xi_2\rangle^{s_2}\widehat{f} \in \ell^2(\widehat{G})\},$$
where $F\in \ell^2(\widehat{G})$ if and only if 
$ \sum_{[\xi]\in \widehat{G} } d_\xi \| F(\xi)\|_{HS}^2<\infty.$
\end{definition}
One may check that the spaces $h^{s_1,s_2}(\widehat{G})$ are indeed complete with respect to the scalar product
$$(f,g)_{s_1,s_2}:=\sum_{[\xi]\in\widehat{G}} d_\xi \langle\xi_1\rangle^{2s_1}\langle\xi_2\rangle^{2s_2}\mathrm{Tr}(\widehat{g}(\xi)^* \widehat{f}(\xi)).$$
Therefore, the Sobolev spaces $H^{s_1,s_2}(G)$ are also complete.

\section{Kernel estimates}
This section is devoted to the proof of some estimates for the (right-convolution) kernels of bisingular pseudodifferential operators on compact Lie groups. 
These estimates will be employed in the next section to develop the global calculus of bisingular operators.

Before proving the estimates, we will first give some properties  representing the suitable bisingular generalization of certain results holding in the standard (global) compact case.
\medskip

\noindent \textbf{Notation.} Recall that $\langle \xi_j\rangle^s:=(1+\lambda_{\xi_j})^{s/2}$, $j=1,2$. Additionally, we assume  $\triangle_1, \triangle_2$ to be the the admissible 
collections of difference operators previously defined. Note that we shall often use the notation $S^{m_1,m_2}(G)$ for  $S^{m_1,m_2}(G\times \widehat{G})$.
\begin{proposition}\label{prop.funct.est}
Let $\triangle^{\alpha,\beta}:=\triangle_1^\alpha\triangle_2^\beta$, then, for any $m_1,m_2\in\mathbb{R}$ and multi-indeces $\alpha\in\mathbb{N}^{n_{P}}$, $\beta\in\mathbb{N}^{n_{R}}$, 
there exists $d\in \mathbb{N}_0$ and $C>0$ such that,
for all $f_1,f_2\in C^d([0,+\infty ))$, $\xi=\xi_1\otimes \xi_2\in\widehat{G}$, and $t_1, t_2\in (0,1)$, we have

$$\| \triangle^{\alpha,\beta}f_1(t_1\lambda_{\xi_1}) f_2(t_2\lambda_{\xi_2})\|_{\mathscr{L}(\mathcal{H}_\xi)}
\leq C t_1^{m_1/2}\langle \xi_1\rangle^{m_1-|\alpha|} \underset{\ell_1=0,...,d}{\sup_{\lambda_{\xi_1}\geq 0}}|\partial^{\ell_1}_{\lambda_{\xi_1}}f_1(\lambda_{\xi_1})|$$
$$ \times 
\,\, t_2^{m_2/2}\langle \xi_2\rangle^{m_2-|\beta|}
\underset{\ell_2=0,...,d}{\sup_{\lambda_{\xi_2}\geq 0}}|\partial^{\ell_2}_{\lambda_{\xi_2}}f_2(\lambda_{\xi_2})|,$$
in the sense that if the supremum on the right-hand side is finite, then the left-hand side is also finite and the inequality holds.
\end{proposition}
\proof [Sketch of proof]
Due to the form of $ \triangle^{\alpha,\beta}$ we have
$$\| \triangle^{\alpha,\beta}f_1(t_1\lambda_{\xi_1}) f_2(t_2\lambda_{\xi_2})\|_{\mathscr{L}(\mathcal{H}_\xi)}=
\| \triangle^{\alpha}_1f_1(t_1\lambda_{\xi_1})\|_{\mathscr{L}(\mathcal{H}_{\xi_1})} 
\| \triangle^{\beta}_2 f_2(t_2\lambda_{\xi_2})\|_{\mathscr{L}(\mathcal{H}_{\xi_2})}.$$
Therefore, by Proposition 6.1 in \cite{Fischer-JFA_2015} applied separately to each term on the right-hand side of the previous identity, the result follows.
\endproof

\begin{lemma}
	Let $k\in\mathcal{D}'(G)$, with $G=G_1\times G_2$ and $n_i=\mathrm{dim}(G_i)$. Then, if $s_1>n_1/2$ and $s_2>n_2/2$,
	$$\|k\|_{L^2(G)}\lesssim \sup_{\xi\in \widehat{G}} \,\, \langle \xi_1\rangle^{s_1/2}  \langle \xi_2\rangle^{s_2/2}\,\, \|\widehat{k}\|_{\mathscr{L}(\mathcal{H}_\xi)}.$$
Hence  $k\in L^2(G)$ when there exist $s_1>n_1/2$ and $s_2>n_2/2$ such that the right-hand side is finite.
\end{lemma}

\proof
Let $B_{s_1,s_2}(x,y)=B_{s_1}(x_1,y_1)\otimes B_{s_2}(x_2,y_2)=B_{s_1}(y_1)\otimes B_{s_2}(y_2)$ as in Lemma \ref{lemma.embedding} ($B_{s_1,s_2}(x,y)$ is independent of $x$). Then, for $s_1,s_2>0$, we can write
$$k(y)=((I_1+L_{G_1})^{s_1/2}\otimes(I_2+L_{G_2})^{s_2/2})(k\ast (B_{s_1}\otimes B_{s_2}) )(y),$$
which gives, in particular, that
$$\widehat{k}(\xi)= \langle \xi_1\rangle^{s_1/2}  \langle \xi_2\rangle^{s_2/2} \widehat{B_{s_1}\otimes B_{s_2}} (\xi) \, \widehat{k}(\xi).$$
Therefore, for $s_1>n_1/2$ and $s_2>n_2/2$, we get
$$\| k\|^2_{L^2(G)} \underset{\text{Plancherel}}{=}\sum_{[\xi]\in\widehat{G}} d_\xi \| \widehat{k}(\xi)\|_{HS}^2$$
$$\leq \sum_{[\xi]\in\widehat{G}} d_\xi \|  \widehat{B_{s_1}\otimes B_{s_2}} (\xi)\|^2_{HS}
\|  \langle \xi_1\rangle^{s_1/2}  \langle \xi_2\rangle^{s_2/2} \widehat{k}(\xi) \|^2_{\mathscr{L}(\mathcal{H}_\xi)}$$
$$\leq \| B_{s_1}\otimes B_{s_2}\|^2_{L^2(G)}\sup_{[\xi]\in\widehat{G}} \langle \xi_1\rangle^{s_1}  \langle \xi_2\rangle^{s_2} \| \widehat{k}(\xi) \|^2_{\mathscr{L}(\mathcal{H}_\xi)}$$
$$\underset{\text{Lemma \ref{lemma.embedding}}}{\lesssim}
\sup_{[\xi]\in\widehat{G}} \langle \xi_1\rangle^{s_1}  \langle \xi_2\rangle^{s_2} \| \widehat{k}(\xi) \|^2_{\mathscr{L}(\mathcal{H}_\xi)},$$
which concludes the proof.
\endproof

\begin{lemma}\label{lemma.kernel1}
	Let $\sigma\in S^{m_1,m_2}(G)$ with (right-convolution) kernel $k_x (\cdot):=k(x,\cdot)$. Then the following properties hold:
	\begin{enumerate}
		\item [1.] The kernel associated with $\partial^{\gamma_1, \gamma_2} \triangle^{\alpha_1,\alpha_2}\sigma \in S^{m_1-|\alpha_1|,m_2-|\alpha_2|}(G)$, for any $\alpha_i\in \mathbb{N}_0^{n_{\triangle_i}}$, and $\gamma_1 \in\mathbb{N}_0^{n_1},\gamma_2\in\mathbb{N}_0^{n_2}$, is given by $q^{\alpha_1,\alpha_2}\partial ^{\gamma_1,\gamma_2}_{x_1,x_2} k_x$;
		\item[2.] If $\sigma_1,\sigma_2$ are two bisingular symbols with kernels $k^1_x$ and $k^2_x$, respectively, then the kernel of the product $\sigma_1 \sigma_2$ is given by $k_x^1\ast k_x^2$.
	\end{enumerate}
\end{lemma}
\proof
The proof of Lemma \ref{lemma.kernel1} follows immediately by the form and the properties of bisingular symbols.
\endproof

As a consequence of Lemma \ref{lemma.kernel1} we get Corollary \ref{cor.kernel.est} below giving a first key estimate for the kernels of  bisingular pseudodifferential operators.
\begin{corollary}\label{cor.kernel.est}
If $\sigma\in S^{m_1,m_2}(G)$, then, for any $\gamma:=(\gamma_1,\gamma_2) \in \mathbb{N}_0^{n_1}\times \mathbb{N}_0^{n_2}$ and 
$\theta=(\theta_1,\theta_2) \in \mathbb{N}_0^{n_1}\times \mathbb{N}_0^{n_2}$ such that, for all $i=1,2$, $\gamma_i+m_i+n_i<|\alpha_i|$,
the function $\partial^{\gamma}_x\partial^{\theta}_z (q^{\alpha_1,\alpha_2}(z)k_x(z))$ is continuous on $G$ and bounded as follows:
$$|\partial^{\gamma}_x\partial^{\theta}_z (q^{\alpha_1,\alpha_2}(z)k_x(z))|\leq C
\sup_{[\xi]\in\widehat{G}}\|\sigma(x,\xi)\|_{S^{m_1,m_2},|\alpha_1|,|\alpha_2|,\gamma},$$
where $\| \cdot \|_{S^{m_1,m_2},|\alpha_1|,|\alpha_2|,\gamma}$ is the suitable seminorm.
The constant $C$ above depends on the parameters $m_i,\triangle, \gamma_i, \vartheta_i$ for all $i=1,2$.
\end{corollary}
\proof
The proof follows from the proof of Corollary 6.5 in \cite{Fischer-JFA_2015} together with  Lemma \ref{lemma.kernel1} and Lemma \ref{lemma.embedding}.
\endproof

Corollary \ref{cor.kernel.est} immediately gives the proposition below.

\begin{proposition}
If $\sigma\in S^{m_1,m_2}(G)$ then the associated kernel $(x,y)\mapsto k_x(y)$ is smooth on $G\times (G\setminus S)$, with $S=\{ x\in G; x_1=e_1\}\cup \{ x\in G; x_2=e_2\}$.
If $\sigma\in S^{-\infty,-\infty}(G)$ is smoothing then the associated kernel $(x,y)\mapsto k_x(y)$ is smooth on $G\times G$. The converse is also true, namely, if $(x,y)\mapsto k_x(y)$ 
is smooth on $G\times G$ then the associated symbol is smoothing, that is, it belongs to $S^{-\infty,-\infty}(G)$.
\end{proposition}

In order to show some estimates for the kernels, we will need to work inside dyadic pieces where the eigenvalues (i.e. the frequencies in this setting) of $L_{G_1}$ and $L_{G_2}$ are localized. 
In that perspective the following lemma will be crucial to understand how the localized symbol and the corresponding kernel behave.

\begin{lemma}\label{lemma.dyadic.behaviour}
Let $\chi\in C_0^\infty(\mathbb{R})$ be a given function with values in $[0,1]$ and $\chi\equiv 1$ in a neighborhood of $0$. Let $\sigma\in S^{m_1,m_2}(G)$ and let $k_x$ be the associated kernel. For  each $\ell_1,\ell_2\in\mathbb{N}$ we define
$$\sigma_{\ell_1,\ell_2}(x,\xi):=\sigma(x,\xi)\chi(\ell_1^{-1}\lambda_{\xi_1})\chi(\ell_2^{-1}\lambda_{\xi_2}).$$
Then $\sigma_{\ell_1,\ell_2}\in S^{-\infty,-\infty}(G)$ and, for any $\gamma=(\gamma_1,\gamma_2)\in \mathbb{N}_0^{n_1}\times \mathbb{N}_0^{n_2}$,
$$\| \sigma_{\ell_1,\ell_2}\|_{S^{m_1,m_2},\gamma}\leq C(G,m_1,m_2,\gamma)\| \sigma\|_{S^{m_1,m_2},\gamma}.$$
Additionally, the kernel $k_{x,\ell_1,\ell_2}(y)$ associated with $\sigma_{\ell_1,\ell_2}$ is smooth on $G\times G$, and, for all $\beta\in \mathbb{N}_0^{n_1+n_2}$, $\partial^{\beta}k_{x,\ell_1,\ell_2}\rightarrow \partial^{\beta}k_{x}$ in $\mathcal{D}'(G)$ uniformly in $x\in G$ as $\ell_1,\ell_2\rightarrow \infty$.
\end{lemma}
\proof
The proof follows the proof of Lemma 6.6 in \cite{Fischer-JFA_2015} with suitable modifications, namely by using the function $\chi(\ell_1^{-1}\lambda_{\xi_1})\chi(\ell_2^{-1}\lambda_{\xi_2})$ as a cutoff function in the proof (note that $(1-\chi(\ell^{-1}\lambda_\pi))$ in \cite{Fischer-JFA_2015} is replaced by $(1-\chi(\ell_1^{-1}\lambda_{\xi_1})\chi(\ell_2^{-1}\lambda_{\xi_2}))$ here), and by replacing the standard Sobolev spaces $H^s(G)$ by the Sobolev spaces $H^{s_1,s_2}(G_1\times G_2)$.
For the sake of completness we shall give the proof of the secod part of the lemma, that is the convergence of the kernels, where a few arrangements are needed.

Let $s_1=\lceil \frac{n_1}{2}\rceil$ and $s_2=\lceil \frac{n_2}{2}\rceil$, where $\lceil\cdot\rceil$ stands for the upper integer part. By using the bisingular Sobolev spaces we get

$$ \|\partial^\beta (k_{x,\ell_1,\ell_2}-k_x) \|_{H^{-s_1-m_1-1,-s_2-m_2-1}}
= \|\partial^\beta (\sigma_{\ell_1,\ell_2}-\sigma) \|_{h^{-s_1-m_1-1,-s_2-m_2-1}}$$
$$= \|(1-\chi(\ell_1^{-1}\lambda_{\xi_1})\chi(\ell_2^{-1}\lambda_{\xi_2}))\partial^\beta \sigma \|_{h^{-s_1-m_1-1,-s_2-m_2-1}}$$
$$\lesssim  \| \langle \xi_1\rangle^{-m_1-1} \langle \xi_2\rangle^{-m_2-1}(1-\chi(\ell_1^{-1}\lambda_{\xi_1})\chi(\ell_2^{-1}\lambda_{\xi_2}))\partial^\beta \sigma \|_{h^{-s_1,-s_2}}$$
$$\lesssim  \| \langle \xi_1\rangle^{-m_1-1} \langle \xi_2\rangle^{-m_2-1}(1-\chi(\ell_1^{-1}\lambda_{\xi_1})\chi(\ell_2^{-1}\lambda_{\xi_2}))\partial^\beta \sigma \|_{L^\infty (\widehat{G})}.$$

Due to the hypothesis on $\chi$, for some $\varepsilon_1, \varepsilon_2>0$, with $0< \varepsilon_1<\varepsilon_2$,  we have that $\chi\equiv 1$ on $[0,\varepsilon_1]$ and $\chi\equiv 0$ on $[\varepsilon_2,+\infty)$. Therefore we get that $(1-\chi(\ell_1^{-1}\lambda_{\xi_1})\chi(\ell_2^{-1}\lambda_{\xi_2}))\not\equiv 0$ in the following three cases
\begin{enumerate}
\item $\lambda_{\xi_1}> \varepsilon_1 \ell_1, \lambda_{\xi_2}> \varepsilon_1 \ell_2$,
\item $\lambda_{\xi_1}> \varepsilon_1 \ell_1, \lambda_{\xi_2}\leq \varepsilon_1 \ell_2$,
\item $\lambda_{\xi_1}\leq \varepsilon_1 \ell_1, \lambda_{\xi_2}> \varepsilon_1 \ell_2$.
\end{enumerate}
Let us start with the proof of the convergence in case $(1)$. 
The inequalities above lead to 
$$\|\partial^\beta (k_{x,\ell_1,\ell_2}-k_x) \|_{H^{-s_1-m_1-1,-s_2-m_2-1}}$$
$$  \leq\underset{ \lambda_{\xi_2}> \varepsilon_\chi \ell_2}{\max_{\lambda_{\xi_1}> \varepsilon_\chi \ell_1}} \|(1-\chi(\ell_1^{-1}\lambda_{\xi_1})\chi(\ell_2^{-1}\lambda_{\xi_2}))\partial^\beta \sigma \|_{h^{-s_1-m_1-1,-s_2-m_2-1}}$$
$$\leq (1+\varepsilon_1 \ell_1)^{-1} (1+\varepsilon_1 \ell_2)^{-1} \| \langle \xi_1\rangle^{-m_1} \langle \xi_2\rangle^{-m_2}\partial^\beta \sigma \|_{L^\infty (\widehat{G})}$$
$$\lesssim (1+\varepsilon_1 \ell_1)^{-1} (1+\varepsilon_1 \ell_2)^{-1}\| \sigma\|_{S^{m_1,m_2},\beta},$$
which gives, in particular, that
$$ \max_{x\in G} \|\partial^\beta (k_{x,\ell_1,\ell_2}-k_x) \|_{H^{-s_1-m_1-1,-s_2-m_2-1}}
\lesssim (1+\varepsilon_1 \ell_1)^{-1} (1+\varepsilon_1 \ell_2)^{-1}\| \sigma\|_{S^{m_1,m_2},\beta}.$$
This finally yields the convergence $\partial^{\beta}k_{x,\ell_1,\ell_2}\overset{\mathscr{D}'}{\rightarrow} \partial^{\beta}k_{x}$ uniformly in $x\in G$ as $\ell_1,\ell_2\rightarrow \infty$. 

For cases $(2)$ and $(3)$ the proof is the same (by reversing the roles of the parameters) and it is similar to the one in the case $(1)$. For completeness we show the steps in case $(2)$, that is, when $\lambda_{\xi_1}> \varepsilon_1 \ell_1$ and $\lambda_{\xi_2}\leq \varepsilon_1 \ell_2$. 
Under these hypotheses we have

$$\|\partial^\beta (k_{x,\ell_1,\ell_2}-k_x) \|_{H^{-s_1-m_1-1,-s_2-m_2-1}}$$
$$  \leq\underset{ \lambda_{\xi_2}\leq \varepsilon_1 \ell_2}{\max_{\lambda_{\xi_1}> \varepsilon_1\ell_1}} \|(1-\chi(\ell_1^{-1}\lambda_{\xi_1})\chi(\ell_2^{-1}\lambda_{\xi_2}))\partial^\beta \sigma \|_{h^{-s_1-m_1-1,-s_2-m_2-1}}$$
$$   \leq\underset{ \lambda_{\xi_2}\leq \varepsilon_1 \ell_2}{\max_{\lambda_{\xi_1}> \varepsilon_1 \ell_1}}  \| \langle \xi_1\rangle^{-m_1-1} \langle \xi_2\rangle^{-m_2-1}(1-\chi(\ell_1^{-1}\lambda_{\xi_1}))\partial^\beta \sigma \|_{h^{-s_1,-s_2}}$$
$$\leq (1+\varepsilon_1 \ell_1)^{-1}  \| \langle \xi_1\rangle^{-m_1} \langle \xi_2\rangle^{-m_2}\partial^\beta \sigma \|_{L^\infty (\widehat{G})}$$
$$\lesssim (1+\varepsilon_1 \ell_1)^{-1} \| \sigma\|_{S^{m_1,m_2},\beta},$$
yielding, as before, the convergence in $\mathcal{D}'$ uniformly in $x$, which completes the proof. 
\endproof

\begin{lemma}\label{lemma.loc-symb}
Let $\sigma\in S^{m_1,m_2}(G)$, and $\eta\in C_0^\infty(\mathbb{R})$. For any $t_1,t_2\in(0,1)$ we define 
$\sigma_{t_1,t_2}(x,\xi):=\sigma(x,\xi)\eta(t_1\lambda_{\xi_1})\eta(t_2\lambda_{\xi_2})$. Then, for any $m_1',m_2'\in\mathbb{R}$ we have
$$\| \sigma_{t_1,t_2}\|_{S^{m'_1,m'_2},\gamma}\leq C t_1^{\frac{m_1-m'_1}{2}}  t_2^{\frac{m_2-m'_2}{2}}\| \sigma\|_{S^{m_1,m_2},\gamma},$$
where $C=C(m_1,m_2,m'_1,m'_2,\gamma,\eta)$ is independent of $\sigma$, $t_1$ and $t_2$.
\end{lemma}

For the proof of Lemma \ref{lemma.loc-symb} see \cite{Fischer-JFA_2015} (Lemma 6.8).
\medskip

We are now ready to prove the main result of this section concerning some estimates for the (right-convolution) kernel of bisingular pseudodifferential operators. 
Let us remark that these estimates are the suitable generalization to our setting of those holding in the standard (non-bisingular) case (see \cite{Fischer-JFA_2015}). 
Note that below we shall denote by $|y|:=d_G(y,e_G)$, where $d_G(\cdot,\cdot)$ is the geodesic distance (and analogously for $|y_j|,$ $j=1,2$). 
Additionally, for any given  $x=(x_1,x_2)\in G$, for a neighborhood of $x$ we shall mean a Cartesian products of the form $U_1\times U_2$, with $U_i$ being a geodesic neigborhood of $x_i$ for $i=1,2$. 



\begin{theorem}\label{prop.kernel.est}
	Let $\sigma\in S^{m_1,m_2}(G)$ and $(x,y)\mapsto k_x(y) \in C^\infty (G\times (G\setminus S))$ be its associated kernel. Then, for $n_i=\mathrm{dim}(G_i), i=1,2$, the following estimates hold
	
	\begin{itemize}
	\item If $n_i+m_i>0$ for $i=1,2$, then there exists $C>0$ and $a,b\in\mathbb{N}$ (independent of $\sigma$) such that for all $y\not\in S$
	
	$$|k_x(y)|\leq C \sup_{\xi\in\widehat{G}}\|\sigma(x,\xi)\|_{S^{m_1,m_2}_{a,b}}  |y_1|^{-n_1-m_1}
	|y_2|^{-n_2-m_2}.$$
	
	\item If $n_i+m_i=0$  for  $i=1,2$, then there exists $C>0$ and $a,b\in\mathbb{N}$ (independent of $\sigma$) such that for all $y\not\in S$
	$$|k_x(y)|\leq C \sup_{\xi\in\widehat{G}}\|\sigma(x,\xi)\|_{S^{m_1,m_2}_{a,b}}  |\ln |y_1||
|\ln	|y_2||.$$
	
	\item If $n_i+m_i<0$ for  $i=1,2$, then $k_x$ is continuous on $G$ and for all $y\not\in S$
$$ |k_x(y)|\leq C \sup_{\xi\in\widehat{G}}\|\sigma(x,\xi)\|_{S^{m_1,m_2}_{0,0}}.$$
	
	\item If $n_i+m_i>0$ and $n_j+m_j=0$ for $i,j\in \{1,2\}, i\neq j$, then
	$C>0$ and $a,b\in\mathbb{N}$ (independent of $\sigma$) such that for all $y\not\in S$
	
		$$|k_x(y)|\leq C \sup_{\xi\in\widehat{G}}\|\sigma(x,\xi)\|_{S^{m_1,m_2}_{a,b}  }
		|y_i|^{-n_i-m_i} |\ln	|y_j||.$$
		
	\item If $n_i+m_i<0$ and $n_j+m_j=0$ for $i,j\in \{1,2\}, i\neq j$, then
$C>0$ and $\gamma_j\in\mathbb{N}^2$ (independent of $\sigma$) being either of the form $\gamma_j=(a_j,0)$ or of the form $\gamma_j=(0,a_j)$, such that for all $y\not\in S$

$$|k_x(y)|\leq C \sup_{\xi\in\widehat{G}}\|\sigma(x,\xi)\|_{S^{m_1,m_2}_{\gamma_j, 0}  }
 |\ln	|y_j||.$$

\item  If $n_i+m_i>0$ and $n_j+m_j<0$ for $i,j\in \{1,2\}, i\neq j$, then
$C>0$ and $\gamma_i\in\mathbb{N}^2$ (independent of $\sigma$, of the same form as above) such that for all $y\not\in S$

	$$|k_x(y)|\leq C \sup_{\xi\in\widehat{G}}\|\sigma(x,\xi)\|_{S^{m_1,m_2}_{\gamma_i,0}  }
|y_i|^{-n_i-m_i} .$$

	\end{itemize}
\end{theorem}

\proof We shall separately analyze the above cases. Let us remark that throughout the proof we shall use the notation $\mathcal{L}_1:=I_1+L_{G_1}$ and $\mathcal{L}_2:=I_2+L_{G_2}$, where $\mathcal{L}_1$ and $\mathcal{L}_2$ are tought of as operators on $G_1$ and $G_2$ respectively, while $L:=\mathcal{L}_1\otimes\mathcal{L}_2$ is defined on $G=G_1\times G_2$.
\medskip

\noindent{\it Case $n_i+m_i>0$}. The estimate in this case trivially follows from Corollary \ref{cor.kernel.est}.
\medskip

\noindent{\it Toolkit}. 
Let $\eta_0,\eta_1\in C_0^\infty(\mathbb{R})$ be supported in $[-1,1]$ and $[1/2,2]$ respectively, taking values in $[0,1]$, and such that
$$\forall \lambda>0 \quad \sum_{\ell=0}^\infty \eta_\ell(\lambda)=1, \quad \text{where}\,\, \eta_\ell(\lambda):=\eta_1(2^{-(\ell-1)}\lambda), \ell\geq 1.$$
Now, for each $\ell_1,\ell_2\in \mathbb{N}_0$, we define
$\sigma_{\ell_1,\ell_2}(x,\xi):=\sigma(x,\xi) \eta_{\ell_1}(\lambda_{\xi_1}) \eta_{\ell_2}(\lambda_{\xi_2})$
(with $\lambda_{\xi_1}$,$\lambda_{\xi_2}$, recall, being the eigenvalues of $\mathcal{L}_1$ and $\mathcal{L}_2$, respectively), and  denote by $k_{x,\ell_1,\ell_2}$ the corresponding kernel.
Notice that, since $\eta_{\ell_1}(\lambda_{\xi_1})\eta_{\ell_2}(\lambda_{\xi_2})$ is smoothing, then $\sigma_{\ell_1,\ell_2}$ is smoothing too. Moreover, also the mapping $(x,y)\mapsto k_{x,\ell_1,\ell_2}(y)=k_x\ast \eta_{\ell_1}(\mathcal{L}_1)\eta_{\ell_2}(\mathcal{L}_2)\delta_{e_1}\otimes\delta_{e_2} $ is smooth, as $(x,y)\mapsto k_{x}(y)$ is smooth on $G\times G\setminus S$ and $\eta_{\ell_1}(\mathcal{L}_1)\eta_{\ell_2}(\mathcal{L}_2)\delta_{e_1}\otimes\delta_{e_2}$  is smooth on $G$.

Observe now that one has the following convergence in $C^\infty(G\times (G \setminus S))$ 
$$k_x(y)=\lim_{N_1,N_2\rightarrow\infty}\sum_{\ell_1=0}^{N_1}\sum_{\ell_2=0}^{N_2}
k_{x,\ell_1,\ell_2}(y)=\Big(k_x\ast \sum_{\ell_1=0}^{N_1}\sum_{\ell_2=0}^{N_2} \eta_{\ell_1}(\mathcal{L}_1)\eta_{\ell_2}(\mathcal{L}_2)\delta_{e_1}\otimes\delta_{e_2}\Big)(y),$$
and that the following bound holds for $y\not\in S$
$$|k_x(y)|\leq \sum_{\ell_1,\ell_2}|k_{x,\ell_1,\ell_2}(y)|.$$

With this in mind we have, by Corollary \ref{cor.kernel.est} and Lemma \ref{lemma.dyadic.behaviour}, that for any given $\alpha_i \in N_0^{n_{\triangle_i}}$, with $i=1,2$, and for any given $m'_i\in\mathbb{R}$, $i=1,2$, such that $m'_i+n_i<|\alpha_i|$,

$$\sup_{x\in G}|q^{\alpha_1,\alpha_2}(z)k_{x,\ell_1,\ell_2}|
\lesssim \sup_{[\xi]\in\widehat{G}} \|\sigma_{\ell_1,\ell_2}(x,\xi)\|_{S^{m'_1,m'_2}_{(|\alpha_1|,|\alpha_2|),0}}$$
\begin{equation}\label{estk}
\lesssim \|\sigma\|_{S^{m_1,m_2}_{(|\alpha_1|,|\alpha_2|),0}}
2^{-(\ell_1-1)\frac{m'_1-m_1}{2}}2^{-(\ell_2-1)\frac{m'_2-m_2}{2}}.
\end{equation}

Note that, for all $z\in G$ and for all $a_1,a_2\in 2\mathbb{N}_0$, we have
$$|z_1|^{a_1}|z_2|^{a_2}\lesssim \sum_{|\alpha_1|=a_1, |\alpha_2|=a_2} |q^{\alpha_1,\alpha_2}(z)|.$$
The previous estimate is of course meaningful in a neghborhhod $U=U_1\times U_2$ where  $U_1$ and $U_2$ are geodesic neigborhoods of $e_1$ and $e_2$ respectively, in which, in the following, we will be working. Note that outside that neighborhood the estimates in the statement are straightforward, because of the smoothness of the kernel.
Therefore, for all $a_1,a_2\in 2\mathbb{N}_0$ and $m'_1,m'_2$ such that $m'_i+n_i<a_i$, $i=1,2$, \eqref{estk} implies
\begin{equation}\label{diagic.estk}
|z_1|^{a_1}|z_2|^{a_2}|k_{x,\ell_1,\ell_2}(z)|\lesssim\|\sigma\|_{S^{m_1,m_2}_{(a_1,a_2),0}}
2^{\ell_1\frac{m_1-m'_1}{2}}2^{\ell_2\frac{m_2-m'_2}{2}}.
\end{equation}

Since we want to study the behaviour of $k_x(y)$ close to the set $S$, we will be considering each of the following situations
\begin{enumerate}
\item[1.] $|z_1|<1$ and $|z_2|<1$;
\item[2.] $|z_1|<1$ and $|z_2|\geq1$ (resp. $|z_1|\geq1$ and $|z_2|<1$).
\end{enumerate}

\medskip

\noindent{\it Case $n_i+m_i>0$ for all $i=1,2$}. 
When $|z_1|<1$ and $|z_2|<1$, we can chose $\ell_{0_i}\in \mathbb{N}_0$ such that
$$2^{-\ell_{0_i}}\leq |z_i| \leq 2^{-\ell_{0_i}+1},\quad  i=1,2.$$
In order to derive the desired estimate we write  
$$\sum_{\ell_1=0}^{N_1}\sum_{\ell_2=0}^{N_2}=\sum_{\substack{\ell_1\leq \ell_{0_1}\\ \ell_2\leq \ell_{0_2} }}+\sum_{\substack{\ell_1\leq \ell_{0_1}\\ N_2\geq \ell_2> \ell_{0_2} }}+\sum_{\substack{N_1\geq \ell_1> \ell_{0_1}\\ \ell_2\leq \ell_{0_2} }}+\sum_{\substack{N_1\geq \ell_1> \ell_{0_1}\\ N_2\geq \ell_2> \ell_{0_2} }},$$
and study the behaviour of $k_{x,\ell_1,\ell_2}$ in the  cases 
\begin{itemize}
\item[1.] $\ell_i\leq \ell_{0_i}$ for $i=1,2$,
\item[2.] $\ell_i> \ell_{0_i}$ for $i=1,2$,
\item[3.]$\ell_1\leq \ell_{0_1}$ and $\ell_2>\ell_{0_2}$ (resp. $\ell_2\leq \ell_{0_2}$ and  $\ell_1>\ell_{0_1}$),
\end{itemize}
separately.

For $\ell_i\leq \ell_{0_i}$, for $i=1,2$,  from \eqref{diagic.estk} we get
$$\sum_{\substack{\ell_1\leq \ell_{0_1}\\ \ell_2\leq \ell_{0_2}}} |k_{x,\ell_1,\ell_2}(z)|\lesssim\|\sigma\|_{S^{m_1,m_2}_{(a_1,a_2),0}}
 |z_1|^{-a_1}2^{\ell_{0_1}\frac{m_1-m'_1}{2}} |z_2|^{-a_2} 2^{\ell_{0_2}\frac{m_2-m'_2}{2}}.$$
We then choose $a_i\in 2\mathbb{N}_0$ and $m'_i\in\mathbb{R}$, for $i=1,2$, such that
\begin{equation}\label{cond.a-m}
m_i+n_i>a_i\geq m_i +n_i-2\quad \text{and}\quad \frac{m_i-m_i'}{2}=m_i+n_i-a_i>0,
\end{equation}
which yields
$$\sum_{\substack{\ell_1\leq \ell_{0_1}\\ \ell_2\leq \ell_{0_2}}} |k_{x,\ell_1,\ell_2}(z)|\lesssim\|\sigma\|_{S^{m_1,m_2}_{(a_1,a_2),0}}
|z_1|^{-a_1-\frac{m_1-m_1'}{2}} |z_2|^{-a_2-\frac{m_2-m_2'}{2}}$$
$$\lesssim \|\sigma\|_{S^{m_1,m_2}_{(a_1,a_2),0}}
|z_1|^{-m_1-n_1} |z_2|^{-m_2-n_2}.$$

For $\ell_i>\ell_{0_i}$ ($\ell_i\leq N_i$) we make a different choice for $a_i$ and $m'_i$ in \eqref{cond.a-m} that we call $a_i', m''_i$ in order to keep the notation $a_i, m_i'$ for the choices we made in the previous case $\ell_i\leq \ell_{0_i}$. We now choose  $a_i'=a_i+2$ and $m''_i$ satisfying 
$$ \frac{m_i-m_i''}{2}=m_i+n_i-a'_i, \quad i=1,2.$$
Since $m_i< m''_i$ now, we have that
$$ \sum_{\substack{N_1\geq \ell_1> \ell_{0_1}\\ N_2\geq \ell_2> \ell_{0_2} }} |k_{x,\ell_1,\ell_2}(z)| \lesssim\|\sigma\|_{S^{m_1,m_2}_{(a_1,a_2),0}}
|z_1|^{-a'_1}2^{\ell_{0_1}\frac{m_1-m_1''}{2}} |z_2|^{-a'_2}2^{\ell_{0_2}\frac{m_2-m_2''}{2}}$$
$$\lesssim \|\sigma\|_{S^{m_1,m_2}_{(a_1,a_2),0}}
|z_1|^{-m_1-n_1} |z_2|^{-m_2-n_2}.$$

For $\ell_1\leq \ell_{0_1}$ and $\ell_2>\ell_{0_2}$ (resp. $\ell_2\leq \ell_{0_2}$ and  $\ell_1>\ell_{0_1}$) we make a different choice of $a_i$ and $m'_i$ that we call $a_i'', m'''_i$ in order to keep the previous notation for the other cases.
By choosing $a_1''=a_1, m_1'''=m_1'$, $a_2''=a_2'$ and $m_2'''=m_2''$ we get, once again from \eqref{diagic.estk}, that
 $$\sum_{\substack{\ell_1\leq \ell_{0_1}\\ N_2\geq \ell_2> \ell_{0_2} }}
|k_{x,\ell_1,\ell_2}(z)|
 \lesssim\|\sigma\|_{S^{m_1,m_2}_{(a_1,a_2),0}}
|z_1|^{-a''_1}2^{\ell_{0_1}\frac{m_1-m_1'''}{2}} |z_2|^{-a''_2}2^{\ell_{0_2}\frac{m_2-m_2'''}{2}}$$
$$\lesssim \|\sigma\|_{S^{m_1,m_2}_{(a_1,a_2),0}}
|z_1|^{-m_1-n_1} |z_2|^{-m_2-n_2}.$$

The estimate in the case when $\ell_2\leq \ell_{0_2}$ and  $\ell_1>\ell_{0_1}$ follows similarly by exchanging the role of $\ell_1$ an $\ell_2$.
Collecting the (four) estimates together we get the desired result (keeping the biggest seminorm) in the case 
when $|z_1|<1$ and $|z_2|<1$.

In the case when $|z_1|<1$ and $|z_2|\geq1$, we can choose $\ell_{0_1}$ as before, and, once again, split the analysis into the cases $\ell_1\leq \ell_{0_1}$ and  $\ell_1>\ell_{0_1}$. Note that we do not split the sum in $\ell_2$ in this case, so we will make a single choice for $a_2$ and $m_2'$. By choosing $a_1$, $m'_1$, $a'_1$, $m''_1$ as before, and $a_2=n_2+m_2+3$, $m_2'=m_2+2$ (so that $m'_2+n_2<a_2$), 
we will get the result in this case (again the result is given in terms of the biggest seminorm).

Finally the case $|z_1|\geq1$ and $|z_2|<1$ is  proved as the last one by reversing the role of $z_1$ and $z_2$. 

Collecting all the estimates above we obtain the result in terms of the biggest seminorm.
\medskip

\noindent{\it Case: $m_i+n_i=0$ for $i=1,2$.}
We consider again all the cases $|z_i|<1$, $|z_j|\geq1$, $i,j=1,2$ ($i\neq j$).
When $|z_1|<1$ and $|z_2|<1$ we fix $\ell_{0_i}$ as before and consider the cases 1, 2 and 3 (and the respective case of the last one) as above.  
Then, for $\ell_i\leq \ell_{0_i}$, for $i=1,2$,  from \eqref{diagic.estk} with $a_i=0$, $m'_i=m_i$, for all $i=1,2$, we get
$$\sum_{\substack{\ell_1\leq \ell_{0_1}\\ \ell_2\leq \ell_{0_2}}} |k_{x,\ell_1,\ell_2}(z)|\lesssim\|\sigma\|_{S^{m_1,m_2}_{a_1,a_2,0}}
\,\, \ell_{0_1}\ell_{0_2}\lesssim\|\sigma\|_{S^{m_1,m_2}_{(a_1,a_2),0}} |\ln |z_1||\,\, |\ln |z_2||.$$
When $\ell_i>\ell_{0_i}$  for $i=1,2$ we choose $a'_i=2$ and $m''_i=m_i-4$ for all $i=1,2,$ and get (from \eqref{diagic.estk} with $a_i', m_i''$)
$$ \sum_{\substack{N_1\geq \ell_1> \ell_{0_1}\\ N_2\geq \ell_2> \ell_{0_2} }} |k_{x,\ell_1,\ell_2}(z)|\lesssim\|\sigma\|_{S^{m_1,m_2}_{(a_1,a_2),0}}.$$
When $\ell_1\leq \ell_{0_1}$ and $\ell_2>\ell_{0_2}$ ($\ell_2\leq \ell_{0_2}$  and $\ell_1>\ell_{0_1}$), by choosing $a_1''=a_1$, $m_1'''=m_1'$ and $a_2''=a_2'$, $m_2'''=m_2''$, we obtain 
$$ \sum_{\substack{\ell_1\leq \ell_{0_1}\\ N_2\geq \ell_2> \ell_{0_2} }} |k_{x,\ell_1,\ell_2}(z)|\lesssim\|\sigma\|_{S^{m_1,m_2}_{(a_1,a_2),0}}|\ln |z_1||.$$
Collecting the estimates together the result when $|z_1|<1$ and $|z_2|<1$ follows.

When  $|z_1|<1$ and $|z_2|\geq1$ we fix again $\ell_{0_1}$ as before. Recall that now we do not split the sum in $\ell_2$ and that we will make a single choice for $a_2$ and $m_2'$ 
in  \eqref{diagic.estk}. Then, using estimate \eqref{diagic.estk} with $a_1$ and $m_1'$ (when $\ell_1\leq \ell_{0_1}$), and $a'_1$ and $m''_1$ (when $\ell_1>\ell_{0_1}$) as in the previous case, 
the result follows by choosing $a_2=n_2+m_2+3=3$ and $m_2'=m_2+2$ (where $m_2+n_2<a_2$ is still satisfied).

The case $|z_1|\geq1$ and $|z_2|<1$ is treated as the previous one reversing the roles of $z_1$ and $z_2$.

Finally, collecting all the cases above, we get the result in terms of the biggest seminorm.
\medskip

\noindent{\it Case $n_i+m_i>0$, $n_j+m_j=0$ for $i,j\in \{1,2\}, i\neq j$}.
To fix ideas suppose $n_1+m_1>0$ and $n_2+m_2=0$ since the other case is treated analogously. We then combine the strategies used in the cases $n_i+m_i>0$ for all $i=1,2$ and $n_i+m_i=0$ for all $i=1,2$.

When $|z_1|<1$ and $|z_2|<1$ we fix again $\ell_{0_i}$ such that $|z_i|\sim 2^{-\ell_{0_i}}$, $i=1,2$. Then, for $\ell_i\leq \ell_{0_i}$, we choose
 $a_i\in 2\mathbb{N}_0$ and $m_i\in\mathbb{R}$, for all $i=1,2$, such that
$$m_1+n_1>a_1\geq m_1 +n_1-2\quad \text{and}\quad \frac{m_1-m_1'}{2}=m_1+n_1-a_1>0,$$
 $$a_2=0,\quad  m'_2=m_2.$$
so that, from \eqref{diagic.estk}, we obtain
 $$\sum_{\substack{\ell_1\leq \ell_{0_1}\\ \ell_2\leq \ell_{0_2}}} |k_{x,\ell_1,\ell_2}(z)|\lesssim \|\sigma\|_{S^{m_1,m_2}_{(a_1,a_2),0}}
 |z_1|^{-m_1-n_1} |\ln |z_2||.$$

For $\ell_i>\ell_{0_i}$, for all $i=1,2$, we apply \eqref{diagic.estk} with $a'_1=a_1+2$, $m''_1$ satisfying the same conditions as $m'_1$ with $a_1'$ in place of $a_1$ (where, recall, $a_1,m'_1$ are the parameters used for $\ell_i\leq \ell_{0_i}$), $a'_2=2$ and $m''_2=m_2-4$. 
We then have
$$ \sum_{\substack{N_1\geq \ell_1> \ell_{0_1}\\ N_2\geq \ell_2>\ell_{0_2}}} |k_{x,\ell_1,\ell_2}(z)| \lesssim\|\sigma\|_{S^{m_1,m_2}_{(a_1,a_2),0}}
|z_1|^{-n_1-m_1}.$$

For $\ell_1\leq \ell_{0_1}$ and $\ell_2>\ell_{0_2}$ ($\ell_2\leq \ell_{0_2}$  and $\ell_1>\ell_{0_1}$) we repeat the strategy used before, that is, we choose $a_1''=a_1$, $m_1'''=m_1'$, $a_2''=a_2'$,  and $m_2'''=m_2''$ in \eqref{diagic.estk} and get 
$$ \sum_{\substack{ \ell_1\leq \ell_{0_1}\\ N_2\geq \ell_2>\ell_{0_2}}} |k_{x,\ell_1,\ell_2}(z)| \lesssim\|\sigma\|_{S^{m_1,m_2}_{(a_1,a_2),0}}
|z_1|^{-n_1-m_1}.$$
Hence, collecting all the estimates we get the result when $|z_1|<1$ and $|z_2|<1$.

When $|z_1|<1$ and $|z_2|\geq1$ the proof follows by considering again only the two cases $\ell_1\leq \ell_{0_1}$ and $\ell_1> \ell_{0_1}$ (here we do not split the sum in $\ell_2$, so $0\leq \ell_2\leq N_2$). Using the same choices as before for $a_1,a'_1, m'_1, m''_1$, and choosing $a_2=3$ and $m'_2=m_2+2$ (so that $m_2'+n_2<a_2$) in \eqref{diagic.estk}, where, recall, $a_1, m_1'$ are the parameters used when $\ell_1\leq \ell_{0_1}$, while $a'_1, m''_1$ are those used for $\ell_1>\ell_{0_1}$ (we make a single choice for $a_2$ and $m'_2$ here), then the desired estimates hold when $|z_1|<1$ and $|z_2|\geq1$.

When $|z_1|\geq1$ and $|z_2|<1$ the result is proved by reversing the roles of $z_1$ and $z_2$ in the last case. 

\medskip
\noindent{\it Cases $n_i+m_i<0$ and $n_j+m_j=0$; $n_i+m_i>0$ and $n_j+m_j<0$ ($i\neq j$)}. 

These cases can be treated as the last one, that is, by combing the strategies used for the other cases in the different regions $|z_i|<1, |z_j|\geq 1$, $i,j=1,2$ ($i\neq j$). The proof is left to the reader.

\endproof
\section{Calculus of bisingular pseudodifferential operators }

In what follows we will use the previous properties to prove a composition formula for bisingular operators.

\begin{theorem}[Composition formula]\label{thm.comp}
	Let $\sigma_A\in S^{m_1,m_2}(G\times \widehat{G})$ and $\sigma_B\in S^{m'_1,m'_2}(G\times \widehat{G})$, and $A:=\mathrm{Op}(a)$ and 
$B=\mathrm{Op}(b)$ the corresponding pseudodifferential operators. Then the symbol $\sigma_{AB}$ of $AB$ is, asymptotically,
	
\begin{equation}\label{asymp-formula}
a\# b(x,\xi):=\sigma_{AB}(x,\xi)\sim \sum_{j\geq 0} c_{m_1+m'_1-j, m_2+m'_2-j}(x,\xi), 
\end{equation}
where  
$$c_{m_1+m'_1-j, m_2+m'_2-j}\in S^{m_1+m'_1-j, m_2+m'_2-j}(G\times \widehat{G}),$$
$$c_{m_1+m'_1-j, m_2+m'_2-j}(x,\xi)=  
d'_{m_1+m'_1-j, m_2+m'_2-j}+d''_{m_1+m'_1-j-1, m_2+m'_2-j}$$
\begin{equation}\label{asympt-terms}
+d'''_{m_1+m'_1-j, m_2+m'_2-j-1},
\end{equation}
$$d'_{m_1+m'_1-j, m_2+m'_2-j}= \sum_{ |\alpha_1|=|\alpha_2|=j}\frac{1}{\alpha_1!\alpha_2!}(\triangle^{\alpha_1,\alpha_2}\sigma_A(x,\xi))\partial^{\alpha_1,\alpha_2}\sigma_B(x,\xi),$$
$$d''_{m_1+m'_1-j-1, m_2+m'_2-j} =\sum_{ |\alpha_2|=j}\frac{1}{\alpha_2!}\left(\triangle^{0,\alpha_2}\sigma_A\circ_{\xi_1}\partial^{0,\alpha_2}
\sigma_B- \sum_{ |\alpha_1|\leq j}\frac{1}{\alpha_1!}(\triangle^{\alpha_1,\alpha_2}\sigma_A(x,\xi))\partial^{\alpha_1,\alpha_2}\sigma_B(x,\xi)\right),$$
and
$$d'''_{m_1+m'_1-j, m_2+m'_2-j-1}
=\sum_{ |\alpha_1|=j}\frac{1}{\alpha_1!}\left(\triangle^{\alpha_1,0}\sigma_A\circ_{\xi_2}\partial^{\alpha_1,0}
\sigma_B - \sum_{ |\alpha_2|\leq j}\frac{1}{\alpha_2!}(\triangle^{\alpha_1,\alpha_2}\sigma_A(x,\xi))\partial^{\alpha_1,\alpha_2}\sigma_B(x,\xi)\right)$$
are such that they belong to $S^{m_1+m'_1-j,m_2+m'_2-j}(G\times \widehat{G})$. In particular the asymptotic formula \eqref{asymp-formula} means that, for any given $N>0$,
$$r_N= \sigma_{AB}-\sum_{j<N}  c_{m_1+m'_1-j, m_2+m'_2-j} \in S^{m_1+m'_1-N,m_2+m'_2-N}(G\times \widehat{G}).$$
\end{theorem}

\proof
Let $A$ and $B$ the operators above, then, by \eqref{Kernel} we have

$$ABf(x)=\int_G(Bf)(xz)R_A(x,z^{-1})dz$$
$$=\int_G f(xy^{-1})\left(\int_G R_B(xz,yz)R_A(x,z^{-1}) dz\right)dy$$
$$\underset{y\rightarrow y^{-1}x}{=}\int_G f(y)\left( \int_GR_B(xz,y^{-1}xz)R_A(x,z^{-1}) dz \right)dy$$
$$=\int_G f(y) R_{AB}(x,y^{-1}x)dy,$$
where
$$R_{AB}(x,y):=\int_G R_B(xz,yz)R_A(x,z^{-1})dz.$$
Since $\sigma_{AB}(x,\xi)=\widehat{R_{A,B}}(x,\xi)$ we have
$$\sigma_{AB}(x,\xi)=\int_G \int_G R_A(x,z^{-1}) R_B(xz,yz)\xi^*(y)  dzdy$$
\begin{equation}\label{comp1}
=\int_G \int_G  R_A(x,z^{-1}) \xi^*(z^{-1})   R_B(xz,yz) \xi^*(yz)dzdy.
\end{equation}
Then we write $R_B(xz,yz)=R_B(x_1z_1, x_2 z_2,y_1z_1, y_2z_2)$  and take the Taylor expansion of $R_B$ with respect to the first variable at $z_1=e_1$, that is,
$$R_B(xz,yz)=\sum_{ |\alpha_1|<N}\frac{1}{\alpha_1!}q^{\alpha_1,0}(z_1^{-1},x_2z_2)\partial^{\alpha_1,0}R_B(x_1, x_2z_2,yz)$$
$$+\sum_{|\alpha_1|=N}
\frac{1}{\alpha_1!}q^{\alpha_1,0}(z_1^{-1},x_2z_2)(R_B)_{\alpha_1}(x_1z_1,x_2z_2,yz),$$
where $q^{\alpha_1,0}(x)=r^{\alpha_1}(x_1)$ is constant with tespect to $x_2$.
Now, taking into account that $q^{\alpha_1,0}(x_1,x_2)$ does not depend on the choice of the second variable and  that $q^{0,\alpha_2}(x_1,x_2)$ does not depend on the choice of the first variable,  we expand the previous quantity with respect to the second variable at $z_2=e_2$ and have
$$R_B(xz,yz)=\sum_{ |\alpha_2|<N}\sum_{ |\alpha_1|<N}\frac{1}{\alpha_1!\alpha_2!}q^{\alpha_1,\alpha_2}(z_1^{-1},z_2^{-1})\partial^{\alpha_1,\alpha_2}R_B(x_1, x_2,yz)$$
$$+\sum_{|\alpha_2|=N}\sum_{ |\alpha_1|<N}
\frac{1}{\alpha_1!\alpha_2!}q^{\alpha_1,\alpha_2}(z_1^{-1},z_2^{-1})(\partial^{\alpha_1,0}R_B)_{\alpha_2}(x_1,x_2,yz),$$
$$+\sum_{ |\alpha_2|<N}\sum_{|\alpha_1|=N}
\frac{1}{\alpha_1!\alpha_2!}q^{\alpha_1,\alpha_2}(z_1^{-1},z_2^{-1})\partial^{0,\alpha_2}(R_B)_{\alpha_1}(x_1,x_2,yz)$$
$$+\sum_{ |\alpha_2|=N}\sum_{|\alpha_1|=N}
\frac{1}{\alpha_1!\alpha_2!}q^{\alpha_1,\alpha_2}(z_1^{-1},z_2^{-1})(R_B)_{\alpha_1,\alpha_2}(x_1z_1,x_2z_2,yz).$$
Therefore we have
\begin{eqnarray*}
\sigma_{AB}(x,\xi)&=&\sum_{ |\alpha_1|<N,|\alpha_2|<N}\frac{1}{\alpha_1!\alpha_2!}\int_{G\times G} \xi^*(z^{-1})q^{\alpha_1,\alpha_2}(z^{-1})R_A(x,z^{-1}) \xi^*(yz)\partial^{\alpha_1,\alpha_2}R_B(x,yz)dzdy\\
&+&\sum_{ |\alpha_1|<N}\frac{1}{\alpha_1!}\int_{G\times G} \Bigg(q^{\alpha_1,0}(z^{-1})R_A(x,z^{-1})\xi^*(z^{-1}) \partial^{\alpha_1,0}R_B(x,yz)\xi^*(yz) \\
&-&\sum_{|\alpha_2|<N}\frac{1}{\alpha_2!}q^{\alpha_1,\alpha_2}(z^{-1})R_A(x,z^{-1})\xi^*(z^{-1}) \partial^{\alpha_1,\alpha_2}R_B(x,yz)\xi^*(yz) \Bigg)dzdy\\
&+&\sum_{ |\alpha_2|<N}\frac{1}{\alpha_2!}\int_{G\times G}\Bigg(q^{0,\alpha_2}(z^{-1})R_A(x,z^{-1})  \xi^*(z^{-1}) \partial^{0,\alpha_2}R_B(x,yz)\xi^*(yz) \\
&-&\sum_{|\alpha_1|<N}\frac{1}{\alpha_1!}q^{\alpha_1,\alpha_2}(z^{-1})R_A(x,z^{-1})\xi^*(z^{-1}) \partial^{\alpha_1,\alpha_2}R_B(x,yz) \xi^*(yz)\Bigg) dzdy\\
&+&\sum_{|\alpha_1|=N,|\alpha_2|=N}\frac{1}{\alpha_1!\alpha_2!}\int_{G\times G} q^{\alpha_1,\alpha_2}(z^{-1})R_A(x,z^{-1}) \xi^*(z^{-1}) (R_B)_{\alpha_1,\alpha_2}(xz,yz)\xi^*(yz)dzdy,
\end{eqnarray*}
and, rearranging the terms, get
$$\sigma_{AB}= \sum_{ |\alpha_1|=|\alpha_2|<N}\frac{1}{\alpha_1!\alpha_2!}(\triangle^{\alpha_1,\alpha_2}\sigma_A(x,\xi)) \partial^{\alpha_1,\alpha_2}\sigma_B(x,\xi)$$
\begin{eqnarray*}
&+&\sum_{ |\alpha_1|<N}\frac{1}{\alpha_1!}\left( (\triangle^{\alpha_1,0}\sigma_A \circ_{\xi_2} \partial^{\alpha_1,0}\sigma_B)(x,\xi) -\sum_{ |\alpha_2|\leq|\alpha_1|}\frac{1}{\alpha_2!} (\triangle^{\alpha_1,\alpha_2}\sigma_A(x,\xi)) \partial^{\alpha_1,\alpha_2}\sigma_B(x,\xi)\right)\\
&+&\sum_{ |\alpha_2|<N}\frac{1}{\alpha_2!}\left( (\triangle^{0,\alpha_2}\sigma_A \circ_{\xi_1} \partial^{0,\alpha_2}\sigma_B)(x,\xi) -\sum_{ |\alpha_1|\leq|\alpha_2|}\frac{1}{\alpha_1!} (\triangle^{\alpha_1,\alpha_2}\sigma_A(x,\xi)) \partial^{\alpha_1,\alpha_2}\sigma_B(x,\xi)\right)\\
&+&\sum_{|\alpha_1|=N,|\alpha_2|=N}\frac{1}{\alpha_1!\alpha_2!}\int_{G\times G}q^{\alpha_1,\alpha_2}(z^{-1})R_A(x,z^{-1}) \xi^*(z^{-1}) (R_B)_{\alpha_1,\alpha_2}(xz,yz)\xi^*(yz) dzdy\\
&=&\sum_{j<N}\left(d'_{m_1+m'_1-j, m_2+m'_2-j}+d'''_{m_1+m'_1-j, m_2+m'_2-j-1}+d''_{m_1+m'_1-j-1, m_2+m'_2-j}\right)
+r_N.
\end{eqnarray*}

In order to complete the proof we only need to show that $r_N\in S^{m_1+m_1'-N,m_2+m_2'-N}(G\times \widehat{G})$ for all $N\in \mathbb{N}_0$, that is, we have to check that
\begin{equation}
\label{est.remainder}
\sup_{x\in G}\|\partial^{\gamma_1,\gamma_2}\Delta^{\beta_1,\beta_2 }r_N(x,\xi)\|_{\mathscr{L}(\mathcal{H}_\xi)}\lesssim \langle\xi_2\rangle^{m_1+m'_1-|\beta_1|-N} \langle\xi_2\rangle^{m_2+m'_2-|\beta_2|-N},
\end{equation}
for all $\gamma_1,\gamma_2, \beta_1,\beta_2$. 
For simplicity we consider the case $\alpha_1=\alpha_2=\beta_1=\beta_2=0$ since the general case follows similarly.
We then write
 $\xi^*(z)=\langle\xi_1\rangle^{-s_1}\langle\xi_2\rangle^{-s_2}(I_1+L_{G_1})_{z_1}^{s_1}\otimes (I_2+L_{G_2})_{z_2}^{s_2}\, \xi^*(z)$, with integers $s_1,s_2\geq 1$, and have, 
after integrating by parts and using the fact that $(R_B)_{\alpha_1,\alpha_2}(x,y)$ is the kernel of a symbol in $S^{m'_1,m'_2}(G\times \widehat{G})$ ,
 \begin{eqnarray*}
r_N(x,\xi)&=&\langle\xi_1\rangle^{-s_1}\langle\xi_2\rangle^{-s_2}\sum_{\substack{|\alpha_1|=N,|\alpha_2|=N \\ |\gamma_1|+|\gamma_2|=2s_1\\ |\tau_1|+|\tau_2|=2s_2}}
c_{\gamma_1,\gamma_2,\tau_1,\tau_2}\frac{1}{\alpha_1!\alpha_2!}\int_{G} \left(\partial _{z}^{\gamma_1,\tau_1}(q^{\alpha_1,\alpha_2}(z^{-1})R_A(x,z^{-1}))\right) \xi^*(z^{-1})\\ &\times& \partial_{z}^{\gamma_2,\tau_2} \widehat{({R}_B)}_{\alpha_1,\alpha_2}(xz,\xi) dz\\
&=&\langle\xi_1\rangle^{-s_1}\langle\xi_2\rangle^{-s_2}\sum_{\substack{|\alpha_1|=N,|\alpha_2|=N \\ |\gamma_1|+|\gamma_2|=2s_1\\ |\tau_1|+|\tau_2|=2s_2}}
c_{\gamma_1,\gamma_2,\tau_1,\tau_2}\frac{1}{\alpha_1!\alpha_2!}\int_{G} \left(\tilde{\partial}_{z^{-1}}^{\gamma_1,\tau_1}q^{\alpha_1,\alpha_2}(z^{-1})R_A(x,z^{-1})\right) \xi^*(z^{-1})\\ &\times& \partial_{z_1=xz}^{\gamma_2,\tau_2} (\widehat{R}_B)_{\alpha_1,\alpha_2}(z_1,\xi) dz,\\
\end{eqnarray*}
where in the second equality we applied the relation beetween left-invariant and right-invariant vector fields given by 
$\partial^{\alpha,\beta}\{\phi(\cdot ^{-1})\}(x)=(-1)^{|\alpha|+|\beta|}(\tilde{\partial}^{\alpha,\beta}\phi)(x^{-1})$ ($\tilde{\partial}$ denoting the right invariant vector field in our notation), and 
used the left invariance of $\partial^{\gamma_2,\tau_2}$.

The previous computations, in particular, give
 \begin{eqnarray*}
\|r_N(x,\xi)\|_{\mathscr{L}(\mathcal{H}_\xi)}&\leq& 
C_{s_1,s_2}  \sum_{\substack{|\alpha_1|=N,|\alpha_2|=N \\ |\gamma_1|+|\gamma_2|=2s_1\\ |\tau_1|+|\tau_2|=2s_2}}
\langle\xi_1\rangle^{-s_1}\langle\xi_2\rangle^{-s_2}
\frac{1}{\alpha_1!\alpha_2!} \int_{G} |\tilde{\partial}_{z^{-1}}^{\gamma_1,\tau_1}q^{\alpha_1,\alpha_2}(z^{-1})R_A(x,z^{-1})|dz\\
&\times&\sup_{z_1\in G} \| \partial^{\gamma_2,\tau_2}_{z_1} (\widehat{R}_B)_{\alpha_1,\alpha_2}(z_1,\xi)\|_{\mathscr{L}(\mathcal{H}_\xi)}\\
&\leq& C_{s_1,s_2} \sum_{\substack{|\alpha_1|=N,|\alpha_2|=N \\ |\gamma_1|+|\gamma_2|=2s_1\\ |\tau_1|+|\tau_2|=2s_2}} \langle\xi_1\rangle^{m'_1-s_1}\langle\xi_2\rangle^{m'_2-s_2}
\frac{1}{\alpha_1!\alpha_2!} \int_{G} |\tilde{\partial}_{z^{-1}}^{\gamma_1,\tau_1}q^{\alpha_1,\alpha_2}(z^{-1})R_A(x,z^{-1})|dz\\
&\times & \|(\widehat{R}_B)_{\alpha_1,\alpha_2}\|_{S^{m'_1,m'_2}_{(2s_1,2s_2)}}.
\end{eqnarray*}
We now assume that $N$ is sufficiently large, namely $N>N_0:=\max\{m_1,m_2\}$, and choose $s_1=N-m_1$ and $s_2=N-m_2$. In this case, by using Proposition \ref{prop.kernel.est}, we obtain 
$$\int_{G} |\tilde{\partial}_{z^{-1}}^{\gamma_1,\tau_1}q^{\alpha_1,\alpha_2}(z^{-1})R_A(x,z^{-1})|dz
\lesssim  \|\tilde{\partial}^{\gamma_1,\tau_1}\Delta^{\alpha_1,\alpha_2} \sigma_A \|_{S^{m_1-N,m_2-N}}\leq
 \|\sigma_A \|_{S^{m_1-N,m_2-N}_{(N,N),(2(N-m_1),2(N-m_2))}},$$
 and, consequently,
 $$ \|r_N(x,\xi)\|_{\mathscr{L}(\mathcal{H}_\xi)\lesssim }\langle\xi_1\rangle^{m_1+m'_1-N}\langle\xi_2\rangle^{m_2-m'_2-N},\quad \forall N>N_0,$$
 which proves \eqref{est.remainder} for every $N>N_0$ when $\gamma_1=\gamma_2=\beta_1=\beta_2=0$. By using similar arguments together with the Leibniz formula one  proves \eqref{est.remainder} in the general form (possibly with a different $N_0$), which, in particular, gives that $r_N\in S^{m_1+m'_1-N,m_2+m'_2-N}$ for every $N>N_0$.
 
 We are now left with proving that $r_N\in S^{m_1+m'_1-N,m_2+m'_2-N}$ for every $N\leq N_0$. Observe that
 \begin{eqnarray*}
 r_N(x,\xi)&=& \sigma_{AB}(x,\xi)- \sum_{j< N} c_{m_1+m'_1-j, m_2+m'_2-j}(x,\xi)\\
 &=& \sigma_{AB}(x,\xi)- \sum_{j< N_0+1} c_{m_1+m'_1-j, m_2+m'_2-j}(x,\xi)
 +\sum_{N\leq j< N_0+1} c_{m_1+m'_1-j, m_2+m'_2-j}(x,\xi)\\
 &=&r_{N_0+1}(x,\xi)+\sum_{N\leq j< N_0+1} c_{m_1+m'_1-j, m_2+m'_2-j}(x,\xi),
  \end{eqnarray*}
therefore, since 
$$r_{N_0+1}\in S^{m_1+m'_1-N_0-1, m_2+m'_2-N_0-1},$$ 
$$\sum_{N\leq j< N_0+1} c_{m_1+m'_1-j, m_2+m'_2-j}\in S^{m_1+m'_1-N, m_2+m'_2-N},$$
and 
$$S^{m_1+m'_1-N_0-1, m_2+m'_2-N_0-1}\subset S^{m_1+m'_1-N, m_2+m'_2-N},$$
we finally get that  $r_N\in S^{m_1+m'_1-N, m_2+m'_2-N}$ for every $N\leq N_0$. This concludes the proof.

\endproof

\begin{theorem}
Let $\sigma\in S^{m_1,m_2}(G\times\widehat{G})$, then the symbol of the operator $\mathrm{Op}(\sigma)^*$, denoted by $\sigma^*$, is asymptotically given by
\begin{equation}\label{asypt.adj}
\sigma^*(x,\xi)\sim \sum_{j\geq 0} c_{m_1-j,m_2-j}(x,\xi),
\end{equation}
where $c_{m_1-j,m_2-j}\in S^{m_1-j,m_2-j}(G\times\widehat{G})$ and
\begin{equation*}
c_{m_1-j,m_2-j}(x,\xi)=d'_{m_1-j,m_2-j}+d''_{m_1-j-1,m_2-j}+d'''_{m_1-j,m_2-j-1},
\end{equation*}
with, using the notations in \eqref{star1} and \eqref{star2} for $\sigma^{*_1}(x,\xi)$ and $\sigma^{*_2}(x,\xi)$,
\begin{eqnarray}
d'_{m_1-j,m_2-j}&=& \sum_{\substack{|\alpha_1|=|\alpha_2|=j}}\frac{1}{\alpha_1!\alpha_2!}
\Delta^{\alpha_1,\alpha_1}\partial^{\alpha_1,\alpha_2}\sigma(x,\xi)^*, \nonumber\\
d''_{m_1-j-1,m_2-j}&=&\sum_{ |\alpha_1|=j}\frac{1}{\alpha_1!}\Big(\Delta^{\alpha_1,0} \partial^{0,\alpha_1}\sigma^{*_1}(x,\xi)
-\sum_{ |\alpha_2|\leq|\alpha_1|}\frac{1}{\alpha_2!}\Delta^{\alpha_1,\alpha_2}\partial^{\alpha_1,\alpha_2}\sigma(x,\xi)^*\Big), \nonumber\\
d'''_{m_1-j,m_2-j-1} &=& \sum_{ |\alpha_2|<j}\frac{1}{\alpha_2!}\Big(\Delta^{0,\alpha_2} \partial^{0,\alpha_2}\sigma^{*_2}(x,\xi)
-\sum_{ |\alpha_1|\leq |\alpha_2|}\frac{1}{\alpha_1!}\Delta^{\alpha_1,\alpha_2}\partial^{\alpha_1,\alpha_2}\sigma(x,\xi)^*\Big), \nonumber\\
\end{eqnarray}
belonging to $S^{m_1-j,n_2-j}(G\times\widehat{G})$. In particular the asymptotic formula \eqref{asypt.adj} means that, for any $N>0$,
$$r_N=\sigma^*-\sum_{j<N}c_{m_1-j,m_2-j}\in S^{m_1-N,m_2-N}(G\times\widehat{G}).$$
\end{theorem}
\proof
The strategy here is similar to the one used for the asymptotic composition formula. 
Notice that, since the kernel of $\sigma^*(x,D)$ satisfies $k_{\sigma^*}(x,v)=\overline{k_\sigma(xv^{-1},v^{-1})}$,  by taking the Fourier transform in the second variable we have
$$\sigma^*(x,\xi)=\int_G \overline{k_\sigma(xv^{-1},v^{-1})} \,\, \xi_1^*(v_1)\otimes\xi^*_2(v_2)dv.$$
We now expand  $\overline{k_\sigma(xv^{-1},v^{-1})}= \overline{k_\sigma(x_1v_1^{-1},x_2v_2^{-1}, v^{-1})}$ in the first variable at $v_1=e_1$, and, afterwards, in the second variable at $v_2=e_2$, and get
\begin{eqnarray*}
	\overline{k_\sigma(xv^{-1},v^{-1})}&=&\sum_{\substack{|\alpha_1|<N, |\alpha_2|<N}}\frac{1}{\alpha_1!\alpha_2!}
	q^{\alpha_1,\alpha_1}(v)\partial^{\alpha_1,\alpha_2}\overline{k_\sigma(x,v^{-1})}\\
	&+&\sum_{\substack{|\alpha_1|<N, |\alpha_2|=N}}	\frac{1}{\alpha_1!\alpha_2!}q^{\alpha_1,\alpha_2}(v)
	\overline{(\partial^{\alpha_1,0}k_\sigma)_{\alpha_2}(x_1,x_2v_2^{-1}, v^{-1})}\\
	&+&\sum_{\substack{|\alpha_1|=N, |\alpha_2|<N}}	\frac{1}{\alpha_1!\alpha_2!}q^{\alpha_1,\alpha_2}(v)
	\overline{(\partial^{0,\alpha_2}(k_\sigma)_{\alpha_1})(x_1v_1^{-1},x_2, v^{-1})}\\
	&+&\sum_{\substack{|\alpha_1|=N, |\alpha_2|=N}}	\frac{1}{\alpha_1!\alpha_2!}q^{\alpha_1,\alpha_2}(v)
	\overline{(k_\sigma)_{\alpha_1,\alpha_2}(x_1v_1^{-1},x_2v_2^{-1}, v^{-1})}.\\
	&=& I+II+III+IV.
\end{eqnarray*}
Now observe that for $II$ we have
\begin{eqnarray}
II&=&\sum_{\substack{|\alpha_1|<N, |\alpha_2|=N}}	\frac{1}{\alpha_1!\alpha_2!}q^{\alpha_1,\alpha_2}(v)
\overline{(\partial^{\alpha_1,0}k_\sigma)_{\alpha_2}(x_1,x_2v_2^{-1}, v^{-1})}\nonumber\\
 &=&\sum_{|\alpha_1|<N}\frac{1}{\alpha_1!}\Big( q^{\alpha_1}(v_1)
\overline{\partial^{\alpha_1,0}k_\sigma(x_1,x_2v_2^{-1}, v^{-1})} -\sum_{ |\alpha_2|<N}\frac{1}{\alpha_2!}q^{\alpha_1,\alpha_2}(v)
\overline{\partial^{\alpha_1,\alpha_2}k_\sigma(x_1,x_2, v^{-1})}\Big),\nonumber
\end{eqnarray}
which shows that $II$ (by the calculus introduced in \cite{MR-VT-Book_2010}) is the kernel of the pseudodifferential operator with symbol
$$ \sum_{ |\alpha_1|<N}\frac{1}{\alpha_1!}\Big(\Delta^{\alpha_1,0} \partial^{\alpha_1,0}\sigma^{*_2}(x,\xi)
-\sum_{ |\alpha_2|<N}\frac{1}{\alpha_2!}\Delta^{\alpha_1,\alpha_2}\partial^{\alpha_1,\alpha_2}\sigma(x,\xi)^*\Big).$$

For the term $III$ with similar arguments one concludes that $III$ is the kernel of 
$$ \sum_{ |\alpha_2|<N}\frac{1}{\alpha_2!}\Big(\Delta^{0,\alpha_2} \partial^{0,\alpha_2}\sigma^{*_1}(x,\xi)
-\sum_{ |\alpha_1|<N}\frac{1}{\alpha_2!}\Delta^{\alpha_1,\alpha_2}\partial^{\alpha_1,\alpha_2}\sigma(x,\xi)^*\Big).$$

For the term $I$ it is immediate to see that it is the kernel of the operator whose symbol is given by
$$\sum_{\substack{|\alpha_1|<N, |\alpha_2|<N}}\frac{1}{\alpha_1!\alpha_2!}
\Delta^{\alpha_1,\alpha_1}\partial^{\alpha_1,\alpha_2}\sigma(x,\xi)^*.$$
Therefore, putting together the properties above and rearranging the terms, one gets 
\begin{eqnarray}
\sigma^*(x,\xi)&\sim& \sum_{\substack{|\alpha_1|=|\alpha_2|<N}}\frac{1}{\alpha_1!\alpha_2!}
\Delta^{\alpha_1,\alpha_1}\partial^{\alpha_1,\alpha_2}\sigma(x,\xi)^*\nonumber\\
&+&\sum_{ |\alpha_1|<N}\frac{1}{\alpha_1!}\Big(\Delta^{\alpha_1,0} \partial^{\alpha_1,0}\sigma^{*_1}(x,\xi)
-\sum_{ |\alpha_2|\leq|\alpha_1|}\frac{1}{\alpha_2!}\Delta^{\alpha_1,\alpha_2}\partial^{\alpha_1,\alpha_2}\sigma(x,\xi)^*\Big)\nonumber\\
&+&\sum_{ |\alpha_2|<N}\frac{1}{\alpha_2!}\Big(\Delta^{0,\alpha_2} \partial^{0,\alpha_2}\sigma^{*_2}(x,\xi)
-\sum_{ |\alpha_1|\leq |\alpha_2|}\frac{1}{\alpha_1!}\Delta^{\alpha_1,\alpha_2}\partial^{\alpha_1,\alpha_2}\sigma(x,\xi)^*\Big)
\nonumber\\
&+& \sum_{\substack{|\alpha_1|=N, |\alpha_2|=N}}\frac{1}{\alpha_1!\alpha_2!}\int q_{\alpha_1,\alpha_2}(v)
\overline{(k_\sigma)_{\alpha_1,\alpha_2}(x_1v_1^{-1},x_2v_2^{-1}, v^{-1})}(\xi_1^*(v_1)\otimes\xi^*_2(v_2))dv\nonumber\\
&=&\sum_{j<N}\left( d'_{m_1-j,m_2-j}+d''_{m_1-j-1,m_2-j}+d'''_{m_1-j,m_2-j-1}\right)+r_N.\nonumber
\end{eqnarray}
In order to complete the proof it remains to show that $r_N\in S^{m_1-N,m_2-N}(G\times\widehat{G})$ which follows by arguments similar to those used in Theorem \ref{thm.comp}. This concludes the proof.
\endproof

\begin{theorem}[Asymptotic expansion]
Let $\sigma_j$ be a sequence of symbols in $S^{m'_j,m''_j}(G\times \widehat{G})$ with $m'_j,m''_j$ decreasing to $-\infty$. Then there exists $\sigma\in S^{m'_0,m''_0}(G\times \widehat{G})$, unique modulo $S^{-\infty,-\infty}$, such that
\begin{equation}
\sigma-\sum_{j=0}^M\sigma_j \in S^{m'_{M+1},m''_{M+1}}(G\times \widehat{G}),\quad \forall M\in\mathbb{N}.
\end{equation}
\end{theorem}

\proof
Let $\psi\in C^\infty(\mathbb{R};[0,1])$ be such that $\psi\equiv 0$ on $(-\infty,1/2)$ and $\psi\equiv 1$ on $(1,\infty)$. Then, by Proposition \ref{prop-product-Leibniz-formula} and Proposition \ref{prop.funct.est}, we have that, for any given $\tilde{m}_1,\tilde{m}_2\in\mathbb{R}$, 

$$\|\Delta^{\alpha,\beta}\partial^{\gamma_1,\gamma_2}\sigma_j(x,\xi)\psi(t_1\lambda_{\xi_1})\psi(t_2\lambda_{\xi_2})\|_{\mathscr{L}(\mathcal{H}_\xi)}\hspace{5cm}$$
\begin{eqnarray}
&\lesssim& \sum_{\substack{|\alpha|\leq |\alpha_1|+|\alpha_2|\leq 2|\alpha|\\ 
|\beta|\leq |\beta_1|+|\beta_2|\leq 2|\beta|}}\|\Delta^{\alpha_1,\beta_1}\partial^{\gamma_1,\gamma_2}\sigma_j(x,\xi)\Delta^{\alpha_2,\beta_2}\psi(t_1\lambda_{\xi_1})\psi(t_2\lambda_{\xi_2})\|_{\mathscr{L}(\mathcal{H}_\xi)}\nonumber\\
&\lesssim& \|\sigma_j\|_{S^{m'_j,m''_j}_{(2\alpha,2\beta), \gamma}}
\sum_{\substack{|\alpha|\leq |\alpha_1|+|\alpha_2|\leq 2|\alpha|\\ 
		|\beta|\leq |\beta_1|+|\beta_2|\leq 2|\beta|}} \langle\xi_1\rangle^{m'_j-|\alpha_1|} \langle\xi_2\rangle^{m''_j-|\alpha_2|}  t_1^{\tilde{m}_1/2}\langle\xi_1\rangle^{\tilde{m}_1-|\alpha_2|} t_2^{\tilde{m}_2/2}\langle \xi_2\rangle^{\tilde{m}_2-|\beta_2|}.\nonumber
\end{eqnarray}
We then choose $\tilde{m}_1=m'_0-m'_j$ and $\tilde{m}_2=m''_0-m''_j$ and get
$$\|\Delta^{\alpha,\beta}\partial^{\gamma_1,\gamma_2}\sigma_j(x,\xi)\psi(t_1\lambda_{\xi_1})\psi(t_2\lambda_{\xi_2})\|_{\mathscr{L}(\mathcal{H}_\xi)}
\lesssim  \|\sigma_j\|_{S^{m'_j,m''_j}_{(|2\alpha|,|2\beta|), (|\gamma_1|,|\gamma_2|)}}
 t_1^{\frac{m'_0-m'_j}{2}}t_2^{\frac{m''_0-m''_j}{2}}\langle \xi_1\rangle^{m'_0-|\alpha|} \langle \xi_2\rangle^{m''_0-|\beta|},
$$
which, in particular, gives that for any given $a=(a_1,a_2)\in \mathbb{N}_0\times  \mathbb{N}_0$ and $b=(b_1,b_2)\in \mathbb{N}_0\times  \mathbb{N}_0$, 
$$\|\sigma_j(x,\xi)\psi(t_1\lambda_{\xi_1})\psi(t_2\lambda_{\xi_2})\|_{S^{m'_0,m''_0}_{a,b}} \leq C_{a,b,m'_0,m''_0,\sigma_j} t_1^{\frac{m'_0-m'_j}{2}}t_2^{\frac{m''_0-m''_j}{2}}.$$
We now choose a decreasing sequence $t_j$ such that
$$t_j\in (0,2^{-j})\quad \text{and} \quad C_{(j,j),(j,j),m'_0,m''_0,\sigma_j} t_j^{\frac{m'_0-m'_j}{2}}t_j^{\frac{m''_0-m''_j}{2}}\leq 2^{-j},$$
and define $\tilde{\sigma}_j(x,\xi):=\sigma_j(x,\xi)\psi(t_j\lambda_{\xi_1})\psi(t_j\lambda_{\xi_2})$. By using the properties above we get, for all $\ell\in\mathbb{N}_0$,
$$\sum_{j=0}^\infty \| \tilde{\sigma}_j\|_{S^{m'_0,m''_0}_{(\ell,\ell),(\ell,\ell)}}
\leq \sum_{j=0}^\ell \| \tilde{\sigma}_j\|_{S^{m'_0,m''_0}_{(\ell,\ell),(\ell,\ell)}}
+\sum_{j=\ell+1}^\infty 2^{-j}<\infty,$$
which implies that $\sigma=\sum_{j=0}^\infty \sigma_j(x,\xi)\psi(t_j\lambda_{\xi_1})\psi(t_j\lambda_{\xi_2}) \in S^{m'_0,m''_0}(G\times \widehat{G})$, and, consequently,  
by taking the sum for $j\geq M$, also that $\sum_{j=M}^\infty \sigma_j(x,\xi)\psi(t_j\lambda_{\xi_1})\psi(t_j\lambda_{\xi_2}) \in S^{m'_M,m''_M}(G\times \widehat{G})$ for all $M\in \mathbb{N}$.
We then have that 
\begin{eqnarray}
\sigma-\sum_{j=0}^{M-1}\sigma_j&=& \sum_{j=0}^\infty \sigma_j(x,\xi)\psi(t_j\lambda_{\xi_1})\psi(t_j\lambda_{\xi_2})-\sum_{j=0}^{M-1}\Big(1-\psi(t_j\lambda_{\xi_1})\psi(t_j\lambda_{\xi_2})+\psi(t_j\lambda_{\xi_1})\psi(t_j\lambda_{\xi_2})\Big)\sigma_j\nonumber\\
&=&-\sum_{j=0}^{M-1}\Big(1-\psi(t_j\lambda_{\xi_1})\psi(t_j\lambda_{\xi_2})\Big)\sigma_j
+\sum_{j=M}^\infty \tilde{\sigma}_j\nonumber\\
\end{eqnarray}
belongs to $S^{m'_M,m''_M}(G\times\widehat{G})$, since,  by Proposition \ref{prop.funct.est},  $1-\psi(t_j\lambda_{\xi_1})\psi(t_j\lambda_{\xi_2})$ is smoothing.
In order to conclude the proof, we just have to show that $\sigma$ is unique up to smoothing operators. 
This last property easily follows by observing that, if $\tau$ is another symbol with the same asymptotic expansion as $\sigma$, then,
 for any given $M\in\mathbb{N}$,
 $$\sigma-\tau=\Big(\sigma-\sum_{j=1}^{M-1}\sigma_j\Big) -\Big(\tau-\sum_{j=1}^{M-1}\sigma_j\Big)\in S^{m'_M,m''_M}(G\times\widehat{G}),$$
which, finally, shows that $\sigma=\tau$ modulo $S^{-\infty,-\infty}(G\times\widehat{G})$ and proves the result.
\endproof

We will now introduce the definition of bielliptic operators and derive, for these objects, the existence of biparametrices.

\begin{definition}\label{def.bielliptic}
Let $a\in S^{m_1,m_2}(G\times \widehat{G})$ and $A=\mathrm{Op}(a)\in L^{m_1,m_2}(G)$. We say that $A$ is bielliptic if
\begin{itemize}
\item[(i)] $a(x,\xi)$ is invertible for all but finitely many $[\xi]\in \widehat{G}$ and, for such $\xi$, its inverse $a(x,\xi)^{-1}$ satisfies 
$$\|a(x,\xi)^{-1}\|_{\mathscr{L}(H_\xi)}\leq \langle \xi_1\rangle^{-m_1}\langle \xi_2\rangle^{-m_2};$$
\item[(ii)] $a(x_1,x_2,D_1,\xi_2)$ is exactly invertible as an operator in $L^{m_1}(G_1)$ for all $(x_2,\xi_2)\in G_2\times \widehat{G}_2$ with inverse in $L^{-m_1}(G_1)$, and, in particular,
$$(a\circ_1 a^{-1})(x_1,x_2, D_1,\xi_2)=\mathrm{Id}_{\mathscr{D}'(G_1)};$$
\item[(iii)] $a(x_1,x_2,\xi_1,D_2)$ is exactly invertible as an operator in $L^{m_2}(G_2)$ for all $(x_2,\xi_2)\in G_2\times \widehat{G}_2$with inverse in $L^{-m_2}(G_2)$, and, in particular, 
$$(a\circ_2 a^{-1})(x_1,x_2,\xi_1, D_2 )=\mathrm{Id}_{\mathscr{D}'(G_2)}.$$
\end{itemize} 
\end{definition}

\begin{theorem}
Let $A\in L^{m_1,m_2}(G)$ be bielliptic. Then there exists $B\in L^{-m_1,-m_2}(G)$ such that 
$$AB=I+K_1,$$
$$BA=I+K_2,$$
where $I:=\mathrm{Id}_{\mathcal{D}'(G)}$ is the identity map and $K_1,K_2$ are smoothing bisingular operators.
\proof
We start with the  proof of the first assertion, namely, the existence of $B$ such that $AB=I+K_1,$ with $K_1$ smoothing.

First observe that, by definition of biellipticity, one has that $a^{-1}\in S^{-m_1,-m_2}(G\times\widehat{G})$. Then, by taking $b_0(x,\xi)=a(x,\xi)^{-1}$ and by using the asymptotic composition formula together with $(ii)$ and $(iii)$ of Definition \ref{def.bielliptic}, we have that $a\# b_0= \boldsymbol{1}-r_{1}$,  with $r_{1}\in S^{-1,-1}(G\times \widehat{G})$ and $\boldsymbol{1}(\xi)=I_{\mathbb{C}^{d_\xi}}$. We now define $b_j:=b_0\# r_j$, with $r_j:=r_1\# r_{j-1}\in S^{-j,-j}(G\times\widehat{G})$ for $j\geq 2$, and have $a\# b_j=( \boldsymbol{1}-r_1)\#r_j$. Then, for $b\sim \sum_{j\geq 0}b_j$, we obtain, for any $k\in\mathbb{N}$,
$$a\#\sum_{j<k} b_j=( \boldsymbol{1}-r_1)\# \Big( \boldsymbol{1}+\sum_{0<j<k} r_j\Big)$$
$$= \boldsymbol{1}+\sum_{0<j<k} r_j-r_1-r_1+r_1\#\sum_{0<j<k} r_j= \boldsymbol{1}-r_k,$$
where, recall, $r_k\in S^{-k,-k}(G\times \widehat{G})$. This, finally, gives that
$$a\# b- \boldsymbol{1}\in S^{-\infty,\infty}(G\times \widehat{G}),$$
which proves the first assertion.

In order to prove the existence of a left parametrix $B$, that is such that $BA=I+K_2$, with $K_2$ smoothing, one proceeds as before, namely, one takes $b_0=a^{-1}$ and defines $b_0\# a- \boldsymbol{1}=-s_1\in S^{-1,-1}(G\times \widehat{G})$ and $s_j:= s_{j-1}\#s_1$ for all $j\geq 2$. Then, taking $b_j:=s_j\#b_0$, the result follows for $b\sim \sum_{j\geq 0} b_j$. This concludes the proof.
\endproof
\end{theorem}

\appendix
\section{}

\begin{lemma}\label{lemma.vabishingOrder}
	Let $G=G_1\times G_2$ be a compact Lie group, with $G_i$, $i=1,2,$ being a compact Lie group of dimension $n_i=\mathrm{dim}(G_i)$, and let also $q\in \mathcal{D}(G)$ and $a_1,a_2\in \mathbb{N}$. Then the fallowing properties are equivalent
	\begin{enumerate}
		\item[1.] For all $(\alpha_1,\alpha_2)\in \mathbb{N}^{n_1}\times\mathbb{N}^{n_2}$, with $|\alpha_i|<a_i$, then $\partial^{\alpha_1,\alpha_2}q(e_G)=0$, that is, $q$ vanishes of order $(a_1-1,a_2-1)$ at $e_G$.
		\item[2.] For any given differential operator $D^{k_1,k_2}:=D_1^{k_1}D_2^{k_2}\in \mathrm{Diff}^{k_1+k_2}(G)$, $D_j^{k_j}\in \mathrm{Diff}^{k_j}(G_j)$, such that $k_i<a_i$, we have $D^{k_1,k_2}q(e_G)=0$.
		\item[3.] There exists a constant $C_q$ such that, for all $x\in G$, we have $q(x)\leq C_q |x_1|^{a_1}|x_2|^{a_2}$.
	\end{enumerate}
\end{lemma}
Lemma \ref{lemma.vabishingOrder}, whose proof is left to the reader, gives a notion of vanishing order of a function suitable in our setting, where, in particular, 
the vanishing order with respect to each variable is considered. For the standard (non adapted to the bisingular case) notion of vanishing orderer of a function see Lemma A.1 in \cite{Fischer-JFA_2015}.

\begin{proposition}\label{prop}
	Let $m_1,m_2\in\mathbb{R}$ and $a_1,a_2\in \mathbb{N}$. For any given function $q\in\mathcal{D}(G)$ vanishing of order $(a_1-1,a_2-1)$ at $e_G$, there exists $d_1,d_2\in\mathbb{N}_0$ such that, for all $f\in C^{d_1}([0,+\infty);C^{d_2}[0,+\infty))$ satisfying
	$$\|f\|_{\mathcal{M}_{m_1,m_2,d_1,d_2}}:=\sup_{\lambda_1,\lambda_2\geq 0, \ell_1=0,...,d_1,\ell_2=0,...,d_2}(1+\lambda_1)^{-m_1+\ell_1}(1+\lambda_2)^{-m_2+\ell_2}|\partial^{\ell_1}_{\lambda_1}\partial^{\ell_2}_{\lambda_2}f(\lambda_1,\lambda_2)|<\infty,$$
	we have
	$$\| \triangle_q{f(t_1\lambda_{\xi_1},t_2\lambda_{\xi_2})}\|_{\mathcal{L}(\mathcal{H}_\pi)}
	\leq C t_1^{m_1/2} t_2^{m_2/2}(1+\lambda_{\xi_1})^{\frac{m_1-a_1}{2}} (1+\lambda_{\xi_2})^{\frac{m_2-a_2}{2}},\quad \forall \pi\in\widehat{G}, \quad t_1,t_2\in(0,1).$$
	The constant $C$ may be chosen as $C'\|f\|_{\mathcal{M}_{m_1,m_2,d_1,d_2}}$, with $C'=C'(m_1,m_2,q,a_1,a_2)$ also depending on the group $G$ but not on $f, t_1,t_2$ and $\xi=\xi_1\otimes\xi_2$.
\end{proposition}
The proof of the proposition is done following that of \cite{Fischer-JFA_2015} and is also left to the reader.

\begin{lemma}\label{lemma.embedding}
	Let $G=G_1\times G_2$ be such that $\mathrm{dim}(G_1)=n_1$. If $s_1>n_1/2, s_2>n_2/2$, then the kernel $\mathcal{B}_{s_1,s_2}$ of the operator $(I_1+L_{G_1})^{-s_1/2}\otimes(I_2+L_{G_2})^{-s_2/2}$ is square integrable and the continuous inclusion $H^{s_1,s_2}(G)\subset C(G)$ holds.
\end{lemma}

\proof[Sketch of the proof of Lemma \ref{lemma.embedding}]

Notice that 

$$B_{s_2,s_2}(x,y)=B_{s_1}(x_1,y_1)\otimes B_{s_2}(x_2,y_2)$$
where
$B_{s_j}(x_j,y_j)$, defined on $G_j\times G_j$, is the kernel of the operator $(I_j+L_{G_j})^{-s_j}$, $j=1,2$.
Then (see Lemma A.5 in \cite{Fischer-JFA_2015}), we have
$$B_{s_j}=\frac{1}{\Gamma(s_j/2)}\int_{t_j=0}^{\infty} t_j^{s_j/2-1}e^{-t_j}p^{(j)}_{t_j}dt_j,$$
where
$$p^{(j)}_{t_j}:=e^{-t_j \Delta_j}\delta_{e_{G_j}},\quad t_j>0,$$
and $\Gamma$ is the gamma function.
Since (see Lemma A.5 in \cite{Fischer-JFA_2015}) for $s_j>n_j/2$
$$\| B_{s_j}\|_{L^2(G_j)}<\infty,\quad j=1,2,$$
we have
$$\|B_s\|_{L^2(G)}=\| B_{s_1}\|_{L^2(G_1)}\| B_{s_2}\|_{L^2(G_2)}<\infty.$$
Finally the Sobolev embedding will follow from the fact that one can write $f$ as 
$$f=\{\left((I_1+L_{G_1})^{-s_1/2}\otimes(I_2+L_{G_2})^{-s_2/2}\right) f\}\ast B_s,$$ for all $f\in H^{s_1,s_2}(G)$ with $s_1>n_1/2$ and $s_2>n_2/2$.
\endproof
	\bibliographystyle{unsrt}
	\bibliography{FedericoParmeggiani-Bisingular-References}
	
\end{document}